\documentclass[12pt,margin=2in]{article}

\usepackage[a4paper,width=150mm,top=25mm,bottom=25mm,bindingoffset=6mm]{geometry}
\usepackage{fancyhdr}
\usepackage{lmodern}
\usepackage[T1]{fontenc}
\usepackage{tikz}\usetikzlibrary{positioning,arrows}
\usetikzlibrary{babel}
\usepackage{graphicx}
\usepackage{caption}
\usepackage{subcaption}
\usepackage{tikz-cd}
\usepackage{amsthm}
\usepackage{amssymb}
\usepackage{amsmath,mathtools}
\usepackage{amsfonts}
\usepackage{enumitem}
\usepackage{bm}
\usepackage{setspace}
\usepackage{mathrsfs}
\usepackage[dvipsnames]{xcolor}
\usepackage[british,UKenglish]{babel}
\usepackage[bookmarks,bookmarksopen,bookmarksnumbered,breaklinks]{hyperref}
\usepackage{authblk}
\usepackage{fontspec}
\usepackage[autostyle]{csquotes}
\usepackage[labelfont=bf,labelsep=space,figurename=Fig.,justification=centering]{caption}

\title{Period collapse of Markov triangles}
\author{Marc Fares\footnote{ORCID: 0009-0006-3976-4747} \\ \textit{Institut de Mathématiques, Université de Neuchâtel\linebreak Rue Emile Argand 11, 2000 Neuchâtel, Switzerland\linebreak \href{mailto:marc.fares@unine.ch}{marc.fares@unine.ch}}}
\date{}

\usepackage{thmtools}
\newtheorem{theorem}{Theorem}[section]
\theoremstyle{plain}
\newtheorem{corollary}[theorem]{Corollary}
\newtheorem{prop}[theorem]{Proposition}

\newtheorem{lemma}[theorem]{Lemma}
\newtheorem*{claim}{Claim}
\newtheorem{assertion}{Assertion}
\theoremstyle{definition}

\newtheorem{fact}[theorem]{Fact}

\declaretheoremstyle[
spaceabove=\baselineskip,
headfont=\normalfont\scshape,
numbered=yes,
bodyfont=\normalfont,
qed={\itshape~$\blacksquare$}
]{exmpstyle}

\declaretheorem[
style=exmpstyle,
title=Definition,
numberwithin=section,
sibling=theorem,
refname={definitionqed},
Refname={definitionqed}
]{definitionqed}

\declaretheoremstyle[
spaceabove=\baselineskip,
headfont=\normalfont\scshape,
numbered=yes,
bodyfont=\normalfont,
qed=
]{exmpstyleo}

\declaretheorem[
style=exmpstyleo,
title=Definition,
numberwithin=section,
sibling=theorem,
refname={definition},
Refname={definition}
]{definition}

\declaretheoremstyle[
spaceabove=\baselineskip,
headfont=\normalfont\bfseries\itshape,
numbered=yes,
bodyfont=\normalfont,
qed={\itshape~$\blacksquare$}
]{exmpstyle2}

\declaretheorem[
style=exmpstyle2,
title=Example,
numberwithin=section,
sibling=theorem,
refname={exampleqed},
Refname={Example}
]{exampleqed}

\declaretheorem[
style=exmpstyle2,
title=Remark,
numberwithin=section,
sibling=theorem,
refname={remarkqed},
Refname={Remark}
]{remarkqed}

\declaretheoremstyle[
spaceabove=\baselineskip,
headfont=\normalfont\bfseries\itshape,
numbered=yes,
bodyfont=\normalfont,
qed=
]{exmpstyleo2}

\declaretheorem[
style=exmpstyleo2,
title=Remark,
numberwithin=section,
sibling=theorem,
refname={remark},
Refname={Remark}
]{remark}

\newcommand{\R}{\mathbb{R}}

\newcommand{\N}{\mathbb{N}}
\newcommand{\Z}{\mathbb{Z}}

\newcommand{\e}{\varepsilon}
\newcommand{\defeq}{\vcentcolon=}
\newcommand{\eqdef}{=\vcentcolon}

\usepackage[backend=biber,maxbibnames=5,maxalphanames=4,style=numeric-comp,url=false,isbn=false,eprint=false]{biblatex}
\DeclareFieldFormat
[article,inbook,incollection,inproceedings,patent,thesis,unpublished]{title}{#1\isdot}
\renewbibmacro{in:}{\ifentrytype{article}{}{\printtext{\bibstring{in}\intitlepunct}}}
\begin{document}
	\maketitle
	\begin{abstract}
		Cristofaro-Gardiner and Kleinman~\cite{gardiner_kleinman} showed the complete period collapse of the Ehrhart quasipolynomial of Fibonacci triangles and their irrational limits, by studying the Fourier--Dedekind sums involved in the Ehrhart function of right-angled rational triangles. We generalize this result using integral affine geometrical methods to all Markov triangles, as defined by Vianna~\cite{vianna}. In particular, we show new occurrences of strong period collapse, namely by constructing for each Markov number~$p$ a two-sided sequence of rational triangles and two irrational limits with quasipolynomial Ehrhart function of period~$p$.
	\end{abstract}
	\paragraph{Acknowledgements.} The problem of studying period collapse phenomena in the Ehrhart theory of Markov triangles was generously proposed by Dan Cristofaro-Gardiner at the Villa Boninchi conference on singular algebraic curves and quanti\-tative symplectic embeddings in September 2025. I thank my doctoral advisor Felix Schlenk for supervising this project as part of my PhD program at the University of Neuchâtel. I also thank Joel Schmitz for his explanations of integral affine geometry, and Nadine Fares and Lama Yaacoub for their helpful comments.
	\newpage
	\section{Introduction and main results}
	\paragraph{Outline of the paper.} We begin by recalling some elements of integral affine geometry in~\S\ref{section_1} and Ehrhart theory in~\S\ref{section_2}, and then define Markov triangles in~\S\ref{section_3}. The main results are Theorems~\ref{limittriangle} and~\ref{limitingbarycentre}, stated in~\S\ref{section_4}. In~\S\ref{section_5} we discuss connections of the Ehrhart theory of our triangles to symplectic embedding problems. All proofs are in~\S\ref{proofs}.
	\subsection{Integral affine geometry}\label{section_1}
	\paragraph{Basic concepts.} Here we define the notions of planar integral affine geometry used in this paper. We refer the interested reader to~\cite{karpenkov2013geometry} for an extensive account of the theory. By a \textit{polygon}~$P$, we mean a two-dimensional convex polytope, i.e.\ the convex hull of finitely many points in~$\R^2$, called the \textit{vertices} of~$P$, while the boundary line segments connecting vertices are called the \textit{edges} of~$P$. A~polygon is said to be \textit{rational} resp.\ \textit{integral} if its vertices have rational resp.\ integral coordinates.

\noindent For rational polygons we can define some additional notions: consider three distinct points~$p_1,p_2,p_3\in\R^2$ with rational coordinates.
\begin{itemize}
	\item The \textit{primitive integral vector}~$v_{1,2}$ along the line segment~$[p_1,p_2]$ is defined to be the unique vector with coprime integer coordinates obtained by multi\-plying~$p_2-p_1$~by~$k\in\R_{>0}$. The \textit{affine length} of~$[p_1,p_2]$ is then~$r\defeq\frac{1}{k}$, namely the unique~$r\in\R_{>0}$ such that~$p_2-p_1=r v_{1,2}$.
	\item The \textit{determinant} of the angle at~$p_1$ is defined as~$|\det(v_{1,2},v_{1,3})|$, namely the absolute value of the determinant of the matrix formed by the primitive integral vectors~$v_{1,2}$ and~$v_{1,3}$ along the segments~$[p_1,p_2]$ and~$[p_1,p_3]$, respectively. The \textit{integral bisector} of the angle at~$p_1$ is the affine line through~$p_1$ spanned by~$v_{1,2}+v_{1,3}$.
\end{itemize}

An affine transformation~$\varphi\in\R^2\rtimes\operatorname{GL}(2,\R)$ is said to be \textit{integral} if it preserves the lattice~$\Z^2$, i.e.~$\varphi(\Z^2)=\Z^2$. The group of integral affine transformations is the semidirect product~$\Z^2\rtimes\operatorname{GL}(2,\Z)$ of the group of translations along integral vectors with the group of integral linear transformations. Two polygons~$P,P'$ are said to be \textit{integrally congruent} if there exists an element~$\varphi\in\Z^2\rtimes\operatorname{GL}(2,\Z)$ taking~$P$ to~$P'$.

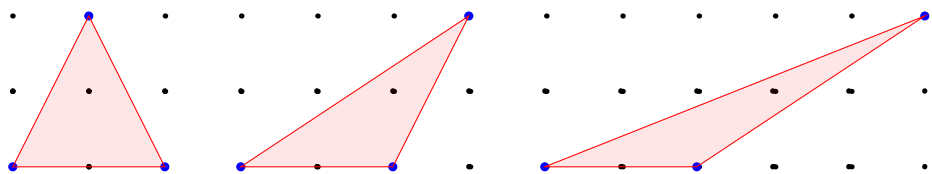
\begin{figure}[h]
	\centering
	\begin{tikzpicture}
		\draw[line width=2pt, line cap=round, dash pattern=on 0pt off 1cm](0,0) grid (12,2);
		\filldraw [blue] (0,0) circle [radius=1.5pt];
		\filldraw [blue] (2,0) circle [radius=1.5pt];
		\filldraw [blue] (1,2) circle [radius=1.5pt];
		\filldraw[color=red,fill opacity=0.1] (0,0)--(2,0)--(1,2)--cycle;
		\filldraw [blue] (3,0) circle [radius=1.5pt];
		\filldraw [blue] (5,0) circle [radius=1.5pt];
		\filldraw [blue] (6,2) circle [radius=1.5pt];
		\filldraw[color=red,fill opacity=0.1] (3,0)--(5,0)--(6,2)--cycle;
		\filldraw [blue] (7,0) circle [radius=1.5pt];
		\filldraw [blue] (9,0) circle [radius=1.5pt];
		\filldraw [blue] (12,2) circle [radius=1.5pt];
		\filldraw[color=red,fill opacity=0.1] (7,0)--(9,0)--(12,2)--cycle;
	\end{tikzpicture}
	\caption{Examples of integrally congruent triangles}
	\label{triangles}
\end{figure}

\paragraph{Half-shears.} Given a vector~$v\in\R^2$, we denote by~$H^+_v$ the half-plane lying to the right of the line spanned by~$v$, and by~$H^-_v$ the one lying to its left; more technically,~$H_v^+$ resp.~$H_v^-$ is the set of vectors~$u\in\R^2$ for which~$\det(u,v)\geq 0$ resp.~$\det(u,v)\leq 0$.
\begin{definitionqed}\label{halfshear}
Given an integral vector~$v=(v_1,v_2)\in\Z^2$, we define the \textit{shear} with respect to~$v$ as integral linear map~$\varphi_v\in\operatorname{GL}(2,\Z)$ given by~$u\mapsto u+\det(v,u)v$, or in matrix form
\[\varphi_v=\begin{pmatrix}
	1-v_1v_2&v_1^2\\
	-v_2^2&1+v_1v_2
\end{pmatrix}.\]
We also define the \textit{half-shear} with respect to~$v$ as the map~$\varphi^{1/2}_v:\R^2\rightarrow\R^2$ with
\[\varphi^{1/2}_v(u)\defeq\begin{cases}
		\varphi_v(u)&\text{for }u\in H_v^+\\
		u&\text{otherwise}
\end{cases}.\]
Note that~$\varphi^{1/2}_v$ is well-defined, since~$\varphi_v$ fixes the line spanned by~$v$.
\end{definitionqed}
In \S\ref{section_3} we will need the following conjugated half-shears.

\begin{definition}\label{geometricmutationdef}
	Let~$P\subset\R^2$ be a rational polygon,~$p\in P$ a vertex of~$P$ and~$v_p$ the primitive integral vector along the integral bisector of the angle at~$p$. The \textit{geometric mutation}~$M_p\defeq M_{p,P}$ of~$P$ at~$p$ is defined as the map obtained by conjugating the half-shear with respect to~$v_p$ by translation along~$p$, i.e.
	\[M_p\defeq T_{p}\circ\varphi^{1/2}_{v_p}\circ T_{-p}.\]
\end{definition}
\begin{remark}\label{shearpreserveslattice}
	The geometric mutation of a rational polygon at an integral vertex preserves the integer lattice. In particular, half-shears preserve the integer lattice.
\end{remark}
	\subsection{Ehrhart theory}\label{section_2}
	\paragraph{Basic notions.} We refer to~\cite[Chapter 3]{beck2020computing} for an excellent exposition of the topic. Let~$P\subset\R^2$ be a polygon.

\begin{definitionqed}
	The \textit{Ehrhart function} of~$P$ is the function in~$t\in\N\cup\{0\}$ defined as the lattice-point enumerator of~$t$-dilates of~$P$; it is denoted by~$L_P$, where
	\[L_P(t)\defeq\#(tP)\cap\Z^2.\qedhere\]
\end{definitionqed}

Two polygons are said to be \textit{Ehrhart equivalent} if they have identical Ehrhart functions. We observe that half-shears preserve Ehrhart functions:

\begin{prop}\label{brunain}
	Let~$v$ be an integral vector, and consider the half-shear~$\varphi^{1/2}_v$ along~$v$. Then~$\varphi^{1/2}_v(P)$ and~$P$ are Ehrhart equivalent.
\end{prop}

The first fundamental result of Ehrhart theory is due to Ehrhart himself~\cite{ehrhart1962polyedres}.
\begin{theorem}[Ehrhart]\label{ehrharttheorem}
	If~$P$ is integral, then~$L_P$ is a polynomial in~$t$.
\end{theorem}

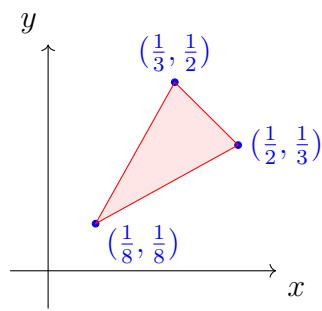
\begin{figure}[h]
	\centering
	\begin{tikzpicture}[scale=5]
		\draw[->] (-0.1,0) -- (0.6,0) node[below right] {$x$};
		\draw[->] (0,-0.1) -- (0,0.6) node[above left] {$y$};
		\filldraw [blue] (1/2,1/3) circle [radius=0.25pt] node[anchor=west] {$(\frac{1}{2},\frac{1}{3})$};
		\filldraw [blue] (1/3,1/2) circle [radius=0.25pt] node[anchor=south] {$(\frac{1}{3},\frac{1}{2})$};
		\filldraw [blue] (1/8,1/8) circle [radius=0.25pt];
		\draw[blue](0.25,0.07)node{$(\frac{1}{8},\frac{1}{8})$};
		\filldraw[color=red,fill opacity=0.1] (1/2,1/3)--(1/3,1/2)--(1/8,1/8)--cycle;
	\end{tikzpicture}
	\caption{Example of a triangle with denominator~$24$}
	\label{denominator}
\end{figure}

If~$P$ is rational, we define its \textit{denominator}~$\operatorname{den}(P)$ to be the least common multi\-ple of the denominators of the coordinates of all vertices of~$P$ (see Figure~\ref{denominator}). \mbox{The~$\operatorname{den}(P)$-dilate} of~$P$ is integral, and hence its Ehrhart function is a polynomial.

\begin{definition}
	A function~$f:\N\cup\{0\}\rightarrow\R$ is said to be \textit{quasipolynomial} if it is periodically polynomial, i.e.\ there exists~$T\in\Z_{>0}$ and a collection of~$T$ poly\-nomials~$\{f_k\}_{k=0}^{T-1}$ such that for all~$t\in\N\cup\{0\}$,~$f(t)=f_k(t)$ whenever~$t\equiv k$~mod~$T$. We call~$T$ the \textit{period} of~$f$.
\end{definition}

\begin{theorem}
	If~$P$ is rational, then~$L_P$ is quasipolynomial with period~$T$ dividing~$\operatorname{den}(P)$.
\end{theorem}

\paragraph{Period collapse.} Given a rational polygon~$P$, we say that \textit{period collapse} occurs if the period~$T$ of~$L_P$ is strictly smaller than~$\operatorname{den}(P)$. Period collapse is said to be \textit{complete} if the period equals~$1$, i.e.~$L_P$ is polynomial, in which case we also say that~$P$~is \textit{pseudo-integral} in accordance with the terminology used by Cristofaro-Gardiner et al.~\cite{cristofarogardiner2015newexamplesperiodcollapse}. Note that while period collapse is not the behaviour that would be expected a priori, it is not a rare phenomenon: McAllister and Woods~\cite{mcallisterwoods} showed that for an arbitrary integer~$D$ and any divisor~$d$ of~$D$, there exists a rational triangle with denominator~$D$ and Ehrhart quasipolynomial with period~$d$. We may also extend the discussion to include irrational polygons; we say that an irrational polygon~$P$ is \textit{pseudo-rational} if~$L_P$ is a quasipolynomial. Note that this is an extreme form of period collapse if we view~$P$ as a rational polygon with infinite denominator.

\begin{exampleqed}\label{fibonacci}
	The Fibonacci sequence~$(\operatorname{Fib}_n)_{n\in\N}$ is defined recursively by
	\[\operatorname{Fib}_1=\operatorname{Fib}_2=1, \, \operatorname{Fib}_{n+2}=\operatorname{Fib}_{n+1}+\operatorname{Fib}_n.\]
	We construct a family of triangles corresponding to the odd-indexed Fibonacci numbers as follows: given a pair~$(r,s)$ of real numbers, the~$(r,s)$-triangle~$\mathcal{F}_{(r,s)}$ is defined by its vertices~$(0,0)$,~$(r,0)$, and~$(0,s)$. The~$n$-th \textit{Fibonacci triangle} is defined as the~$(g_n\defeq \frac{\operatorname{Fib}_{2n-1}}{\operatorname{Fib}_{2n-3}},\frac{\operatorname{Fib}_{2n-3}}{\operatorname{Fib}_{2n-1}})$-triangle, for~$n\geq 2$\footnote{We define the first Fibonacci triangle to be the standard~$2$-dimensional simplex~$\mathcal{F}_{(1,1)}$.}. Recall that the sequence~$\operatorname{Fib}_{n+1}/\operatorname{Fib}_n$ converges to the golden ratio~$\tau=\frac{1+\sqrt{5}}{2}$ as~$n\rightarrow\infty$, which allows us to define the limiting triangle~$\mathcal{F}_\infty\defeq\lim_{n\rightarrow\infty}\mathcal{F}_{(g_n,1/g_n)}=\mathcal{F}_{(\tau^2,1/\tau^2)}$ (see Figure~\ref{limitfibonacci}). $\mathcal{F}_\infty$ is irrational. Cristofaro-Gardiner and Kleinman~\cite{gardiner_kleinman} showed that complete period collapse occurs for all Fibonacci triangles as well as their limit~$\mathcal{F}_\infty$.
\end{exampleqed}
	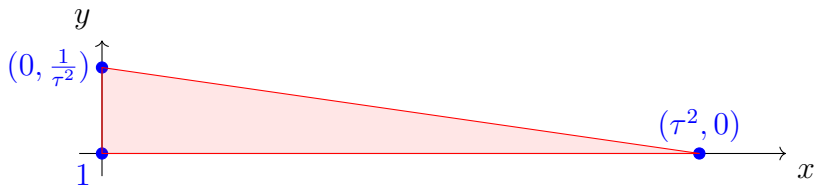
\begin{figure}[h]
		\centering
		\begin{tikzpicture}[scale=3]
			\draw[->] (-0.1,0) -- (3,0) node[below right] {$x$};
			\draw[->] (0,-0.1) -- (0,0.5) node[above left] {$y$};
			\filldraw [blue] (0,0) circle [radius=0.7pt] node[anchor=north east] {$1$};
			\filldraw [blue] (2.62,0) circle [radius=0.7pt] node[anchor=south] {$(\tau^2,0)$};
			\filldraw [blue] (0,0.38) circle [radius=0.7pt] node[anchor=east]{$(0,\frac{1}{\tau^2})$};
			\filldraw[color=red,fill opacity=0.1] (0,0)--(2.62,0)--(0,0.38)--cycle;
		\end{tikzpicture}
		\caption{The limiting Fibonacci triangle}
		\label{limitfibonacci}
	\end{figure}
The goal of this paper is to generalize the period collapse phenomenon for Fibonacci triangles to a larger class of triangles, the so-called \textit{Markov triangles}.
	\subsection{Markov triangles}\label{section_3}
	\paragraph{The Markov tree.} The Diophantine equation
\begin{equation}\label{markov}
p_1^2+p_2^2+p_3^2=3p_1p_2p_3
\end{equation}
is called \textit{Markov's equation}, and a solution~$(p_1,p_2,p_3)\in\Z_{\geq 1}^3$ of \eqref{markov} is called a \textit{Markov triple}, with elements called \textit{Markov numbers}. We will need the following fact from~\cite[Corollary 3.4]{aigner2013markov}.
\begin{fact}
	The elements of a Markov triple are pairwise coprime.
\end{fact}

\begin{figure}[h]
	\centering
\begin{tikzpicture}[level distance=1cm,
	level 3/.style={sibling distance=5cm},
	level 4/.style={sibling distance=2.5cm}]
	\node {$(1,1,1)$}
	child{node{$(1,2,1)$}
		child{node{$\underline{\overline{(1,5,2)}}$}
			child{node{$\overline{(1,13,5)}$}
				child{node{$(1,34,13)$}}
				child{node{$\overline{(5,194,13)}$}}
			}
			child{node{$\underline{(2,29,5)}$}
				child{node{$\underline{(5,433,29)}$}}
				child{node{$(2,169,29)$}}
			}
		}
	};
\end{tikzpicture}
\caption{The first five generations of the Markov tree}
\label{markovtree}
\end{figure}
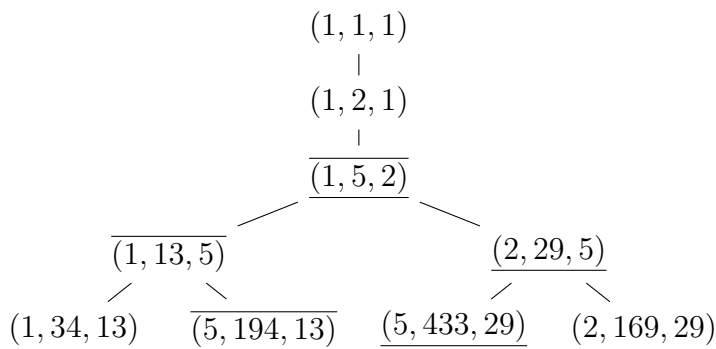

The Markov triples may be generated recursively as follows: replacing the number~$p_2$~by~$\widehat{p_2}\defeq 3p_1p_3-p_2$ in a Markov triple~$(p_1,p_2,p_3)$ yields a new Markov triple. This~$p_2\rightarrow\widehat{p_2}$ trading operation is called \textit{mutation} of~$(p_1,p_2,p_3)$ at~$p_2$. If we assume that~$\max(p_1,p_2)\leq p_3$, then~$\widehat{p_3}\leq p_3\leq\min(\widehat{p_1},\widehat{p_2})$, i.e.\ mutating at~$p_1$ or~$p_2$ increases the largest number in the Markov triple while mutating at~$p_3$ decreases it. We will call a mutation of a Markov triple \textit{increasing} resp.\ \textit{decreasing} if the largest number of the triple increases resp.\ decreases after applying the mutation. We also note that mutating~$(p_1,\widehat{p_2},p_3)$ at~$\widehat{p_2}$ returns the initial triple~$(p_1,p_2,p_3)$. We define a graph with vertices the Markov triples, each two of which are connected by an edge if one can be obtained from the other by a single mutation; our graph is then a trivalent tree, called the \textit{Markov tree}. Plotting the vertices so that an increasing mutation is represented by a downward-pointing edge yields the following representation.\footnote{Note that the three mutations of~$(1,1,1)$ are identical, as well as the two mutations of~$(1,2,1)$ at the~$1$'s; we therefore omit the repeated Markov triples in Figure~\ref{markovtree}.}

Let~$p\notin\{1,2\}$ be a Markov number and~$\mathfrak{m}$ a Markov triple in which~$p$ is the largest number. We denote by~$\Lambda(\mathfrak{m})$ the subgraph of the Markov tree consisting of all triples obtained by an iteration of mutations of~$\mathfrak{m}$ preserving~$p$.~$\Lambda(\mathfrak{m})$ is a bivalent subtree rooted in~$\mathfrak{m}\eqdef(a,p,b)$ consisting of two branches: the~$B_p^a$-branch obtained by mutating~$(a,p,b)$ first at~$b$, and the~$B_p^b$-branch obtained by mutating~$(a,p,b)$ first at~$a$, each of which then extends by mutations that preserve~$p$ and the largest element of the triple. As an example, the subtree formed by the vertices which are either underlined or overlined in Figure~\ref{markovtree} is~$\Lambda(1,5,2)$. The first two vertices of the branch~$B_5^1$ are the overlined triples and those of the branch~$B_5^2$ are the underlined triples. 

In the case~$p=1$, since we omitted two of the three identical subtrees adjacent to~$(1,1,1)$, the subgraph rooted in~$(1,1,1)$ obtained via mutations preserving~$1$ is a branch, with a particularity: given a Markov triple~$(1,p_2,p_3)$, we have~$\widehat{p_2}=3p_3-p_2$, which is the recursive equation defining odd-indexed Fibonacci numbers,\linebreak i.e.~$\operatorname{Fib}_{2n+3}=3\operatorname{Fib}_{2n+1}-\operatorname{Fib}_{2n-1}$.~Since~$\operatorname{Fib}_1=1$~and~$\operatorname{Fib}_3=2$, it follows by induction that the triples of this branch are given by the odd-indexed Fibonacci~numbers, i.e the~$n$-th triple (for~$n\geq 2$) is 
\[(1,\operatorname{Fib}_{2n-1},\operatorname{Fib}_{2n-3}).\]
Accordingly, we call this branch the \textit{Fibonacci branch} of the Markov tree. A similar situation is observed in the case~$p=2$, where the subgraph obtained by mutations preserving~$2$ and rooted in~$(1,2,1)$ is a branch with~$n$-th triple (for~$n\geq 2$)
\[(1,\operatorname{Pell}_{2n-1},\operatorname{Pell}_{2n-3})\]
where~$(\operatorname{Pell}_n)_{n\in\N}$ is the \textit{Pell sequence}, defined recursively by
\[\operatorname{Pell}_0=0,\, \operatorname{Pell}_1=1, \, \operatorname{Pell}_{n+2}=2\operatorname{Pell}_{n+1}+\operatorname{Pell}_n.\]
Accordingly, we call this branch the \textit{Pell branch} of the Markov tree. Note that the Fibonacci and Pell branches are the peripheral branches of the tree. Many more interesting results and open problems on Markov numbers are surveyed in~\cite{aigner2013markov}.

\paragraph{Markov triangles.} Motivated by the complete period collapse of Fibonacci triangles, we study a larger class of rational triangles representing all Markov triples defined by Vianna~\cite{vianna} and Evans~\cite[Appendix I]{evans2023lectures}.
\begin{definitionqed}\label{defmarkov}
	Let~$(p_1,p_2,p_3)$ be a Markov triple. A \nolinebreak{\textit{Markov triangle}~$\Delta_{(p_1,p_2,p_3)}$} representing~$(p_1,p_2,p_3)$ is defined recursively as follows:
	\begin{itemize}
		\item~$\Delta_{(1,1,1)}$ is any triangle with edges having affine length~$1$ and vertices having angles of determinant~$1$.
		\item For the~$n+1$-st generation Markov triple~$(p_1,\widehat{p_2},p_3)$ obtained by mutating the~$n$-th generation triple~$(p_1,p_2,p_3)$ at~$p_2$,~$\Delta_{(p_1,\widehat{p_2},p_3)}$ is defined as the image of the geometric mutation of~$\Delta_{(p_1,p_2,p_3)}$ at the~$p_2$-vertex, namely
		\[\Delta_{(p_1,\widehat{p_2},p_3)}\defeq M_{p_2}(\Delta_{(p_1,p_2,p_3)})\]
		where~$M_{p_2}\defeq M_{p_2,\Delta_{(p_1,p_2,p_3)}}$ is given in Definition~\ref{geometricmutationdef}.\qedhere
	\end{itemize}
\end{definitionqed}

\begin{figure}[h]
	\centering
	\begin{tikzpicture}
		\filldraw [blue] (0,0) circle [radius=1.5pt] node[anchor=north] {$2$};
		\filldraw [blue] (7.5,0) circle [radius=1.5pt] node[anchor=north] {$1$};
		\filldraw [blue] (2.7,1.2) circle [radius=1.5pt] node[anchor=south]{$5$};
		\filldraw[color=red,fill opacity=0.1] (0,0)--(7.5,0)--(2.7,1.2)--cycle;
		\draw[color=blue,dotted] (7.5,0)--(135/58,30/29);
		
		\draw[blue] (3.75,-0.75) node[anchor=west] {$M_1$};
		\draw [blue,->](3.75,-0.25) -- (3.75,-1.25);
		
		\filldraw [blue] (0,-2.5) circle [radius=1.5pt] node[anchor=north] {$2$};
		\filldraw [blue] (8.7,-2.5) circle [radius=1.5pt] node[anchor=north] {$5$};
		\filldraw [blue] (135/58,30/29-2.5) circle [radius=1.5pt] node[anchor=south]{$29$};
		\filldraw[color=red,fill opacity=0.1] (0,-2.5)--(8.7,-2.5)--(135/58,30/29-2.5)--cycle;
		\draw[color=blue,dotted] (7.5,-2.5)--(135/58,30/29-2.5);
	\end{tikzpicture}
	\caption{Geometric mutation of a~$(1,5,2)$-triangle at the~$1$-vertex}
	\label{mutation}
\end{figure}
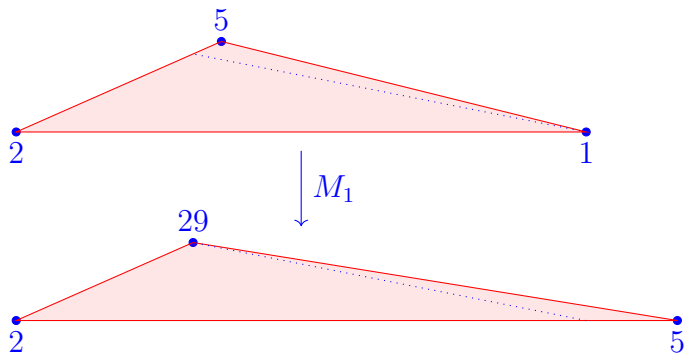

We say that two triangles~$\Delta$ and~$\Delta'$ are \textit{Markov equivalent} if they represent the same Markov triple.

\begin{remark}
	It is not true in general that the image of a triangle~$\Delta$ under geometric mutation~$M_p$ at a vertex~$p$ is again a triangle: it may be a quadrilateral. However in the case of Markov triangles,~$M_p(\Delta)$ is a triangle, as shown by Evans~\cite[Theorem 8.21]{evans2023lectures}. We also give a computational proof of this fact in~\S\ref{section_23}.
\end{remark}
\begin{remarkqed}\label{congruentmarkov}
	The equivalence classes of Markov triangles are invariant\linebreak under the~$\R^2\rtimes\operatorname{GL}(2,\Z)$-action, on the one hand because geometric mutations commute with the~$\operatorname{GL}(2,\Z)$-action, and on the other hand because the class of Markov triangles representing~$(1,1,1)$ is by definition invariant under~$\R^2$-translations.
\end{remarkqed}

It will be useful for our purposes to fix a model triangle in each Markov class.
\begin{definition}
	Let~$\mathfrak{m}\defeq(p_1,p_2,p_3)$ be a Markov triple. The \textit{companion numbers}~$q_1^{\pm}\in\Z_{p_1}$ of~$p_1$ with respect to~$\mathfrak{m}$ are~$q_1^{\pm}\vcentcolon\equiv\pm 3p_3p_2^{-1}$~mod~$p_1$.
\end{definition}
\begin{remarkqed}\label{ordercompanion}
	The set~$\{q_1^\pm\}$ does not depend on the order of the triple~$\mathfrak{m}$, since the congruence~$p_2^2+p_3^2\equiv 0$~mod~$p_1$ 
	implies that~$3p_2p_3^{-1}\equiv -3p_3p_2^{-1}$~mod~$p_1$.
\end{remarkqed}

We will show in Proposition~\ref{pseudomarkov} that every Markov triangle~$\Delta$ representing a triple~$(p_1,p_2,p_3)$ is~$\R^2\rtimes\operatorname{GL}(2,\Z)$-congruent to the triangle~$\Delta(q_1)$ with~$p_1$-,~$p_3$-, and~$p_2$-vertices
\begin{equation}\label{aarous}
	\Bigg\{(0,0),\Big(\frac{(p_1q_1-1)p_2}{p_1p_3},\frac{p_1p_2}{p_3}\Big),\Big(\frac{p_3}{p_1p_2},0\Big)\Bigg\}
\end{equation}
for any companion number~$q_1\vcentcolon\equiv q_1^+\equiv 3p_3p_2^{-1}$~mod~$p_1$ of~$p_1$.
\begin{definitionqed}\label{standardp}
	A Markov triangle representing the triple~$(p_1,p_2,p_3)$ with co\-ordinates given by \eqref{aarous} is said to be in a \textit{standard~$p_1$-position}.
\end{definitionqed}
This definition allows us to view the Fibonacci triangles from Example~\ref{fibonacci} as Markov triangles representing the triples~$(1,\operatorname{Fib}_{2n-1},\operatorname{Fib}_{2n+1})$ in standard~$1$-position.

\begin{remark}\label{equimarkov}
	The~$\operatorname{GL}(2,\Z)$-class of~$\Delta(q_1)$ depends only on the congruence class of~$q_1$~modulo~$p_1$: for~$q_1'\in\Z$ with~$q_1-q_1'=kp_1$,~$\Delta(q_1')$ is the image of~$\Delta(q_1)$ under
	\[\begin{pmatrix}
		1&-k\\0&1
	\end{pmatrix}.\]
	Moreover, the~$\operatorname{GL}(2,\Z)$-class of~$\Delta(q_1)$ does not depend on the order of~$(p_1,p_2,p_3)$: the triangle~$\Delta'$ obtained by swapping~$p_2$ and~$p_3$ in the definitions of~$\Delta(q_1)$ and of~$q_1$ is the image of~$\Delta(q_1)$ under
	\nopagebreak
	\[\psi_{p_1}^{p_2,p_3}\defeq\begin{pmatrix}
		-1-q_1p_1&q_1^2\\
		p_1^2&1-q_1p_1
	\end{pmatrix}.\]
\end{remark}
	\subsection{Ehrhart theory of Markov triangles}\label{section_4}
	In the following, we define the \textit{integral barycentre} of a triangle to be the intersection of its integral bisectors, which is well-defined by Lemma~\ref{barycentre}. Recall that the Ehrhart function of a triangle~$\Delta$ does not change under a half-shear of~$\Delta$, but might change under an~$\R^2$-translation. Given a Markov triple~$(p_1,p_2,p_3)$, we consider:
\begin{itemize}
	\item Markov triangles representing~$(p_1,p_2,p_3)$ in a standard~$p_1$-position,
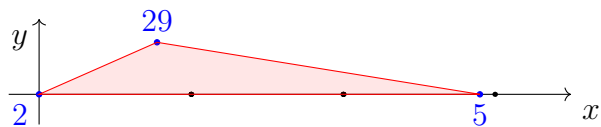
\begin{figure}[h]
	\centering
	\begin{tikzpicture}[scale=2]
		\draw[->] (-0.2,0) -- (3.5,0) node[below right] {$x$};
		\draw[->] (0,-0.2) -- (0,0.5) node[below left] {$y$};
		\draw[line width=2pt, line cap=round, dash pattern=on 0pt off 2cm](0,0) grid (3,0);
		\filldraw [blue] (0,0) circle [radius=0.5pt] node[anchor=north east] {$2$};
		\filldraw [blue] (2.9,0) circle [radius=0.5pt] node[anchor=north] {$5$};
		\filldraw [blue] (45/58,10/29) circle [radius=0.5pt] node[anchor=south]{$29$};
		\filldraw[color=red,fill opacity=0.1] (0,0)--(2.9,0)--(45/58,10/29)--cycle;
	\end{tikzpicture}
	\caption{A~$(2,5,29)$-triangle in standard~$2$-position}
	\label{standardpsetting}
\end{figure}
	
	\item Markov triangles representing~$(p_1,p_2,p_3)$ in a \textit{barycentric position}, i.e.\ with the integral barycentre at the origin of the~$\R^2$-plane. In particular, such a Markov triangle~$\Delta$ obtained by translating a triangle representing~$(p_1,p_2,p_3)$ in standard~$p_1$-position (so that~$\Delta$ has an edge parallel to the~$x$-axis) is said to be in a \textit{standard~$p_1$-barycentric position}.\footnote{The triangles depicted in Figure~\ref{barycentricsetting} are only illustrative; typical Markov triangles in barycentric position are very thin and have one vertex extremely close to the origin, and hence their depiction would not help visualization.}
	\begin{figure}[h]\centering
		\begin{subfigure}[t]{0.5\textwidth}
		\centering
		\begin{tikzpicture}[scale=2]
			\draw[->] (-1.4,0) -- (1.4,0) node[below right] {$x$};
			\draw[->] (0,-1.4) -- (0,1.4) node[below left] {$y$};
			\draw[line width=2pt, line cap=round,dash pattern=on 0pt off 2cm](-1,-1) grid (1.1,1.1);
			\filldraw [blue] (1.1,-1.1) circle [radius=0.5pt] node[anchor=north] {$2$};
			\filldraw [blue] (0.5,0.5) circle [radius=0.5pt] node[anchor=south] {$5$};
			\filldraw [blue] (-1.3,0.1) circle [radius=0.5pt] node[anchor=south east]{$29$};
			\filldraw[color=red,fill opacity=0.1] (1.1,-1.1)--(0.5,0.5)--(-1.3,0.1)--cycle;
		\end{tikzpicture}
		\end{subfigure}
		~
		\begin{subfigure}[t]{0.3\textwidth}
			\centering
			\begin{tikzpicture}[scale=2]
				\draw[->] (-1.4,0) -- (1.4,0) node[below right] {$x$};
				\draw[->] (0,-1.4) -- (0,1.4) node[below left] {$y$};
				\draw[line width=2pt, line cap=round,dash pattern=on 0pt off 2cm](-1,-1) grid (1.1,1.1);
				\filldraw [blue] (1.35,-0.25) circle [radius=0.5pt] node[anchor=north] {$29$};
				\filldraw [blue] (-0.4,1) circle [radius=0.5pt] node[anchor=south] {$5$};
				\filldraw [blue] (-1.3,-0.25) circle [radius=0.5pt] node[anchor=north]{$2$};
				\filldraw[color=red,fill opacity=0.1] (1.35,-0.25)--(-0.4,1)--(-1.3,-0.25)--cycle;
			\end{tikzpicture}
		\end{subfigure}
		\caption{$(2,5,29)$-triangles in a generic barycentric position (left) and in a standard~$2$-barycentric position (right)}
		\label{barycentricsetting}
	\end{figure}
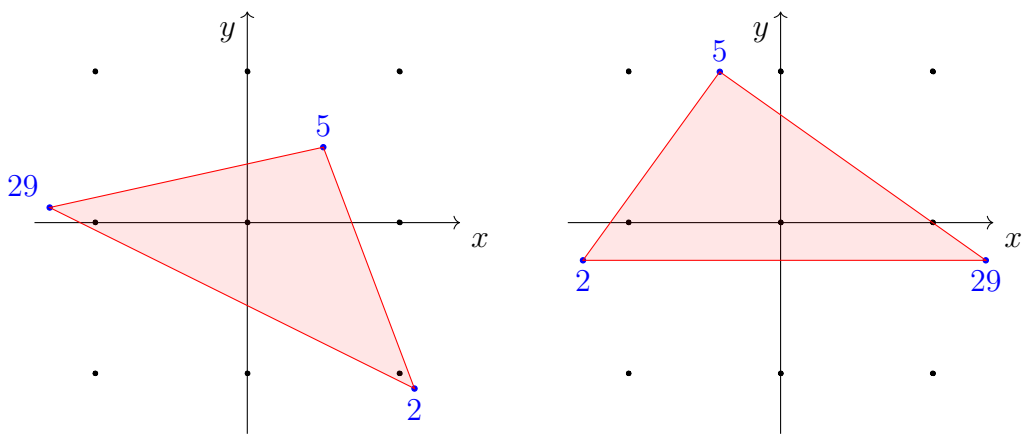
\end{itemize}
\begin{remark}
	Since every~$(p_1,p_2,p_3)$-triangle in barycentric position is \nolinebreak{$\operatorname{GL}(2,\Z)$-congruent} to one in standard~$p_1$-barycentric position, it will suffice for us to consider the latter position.
\end{remark}
\begin{remark}\label{translatebarycentre}
	More precisely, a triangle~$\Delta$ representing the triple~$(p_1,p_2,p_3)$ in standard~$p_1$-barycentric position is obtained from a \mbox{triangle~$\Delta(q_1)$ in standard~$p_1$-position} by translation along the vector pointing from the integral barycentre~$\beta^{(p_1)}$ of~$\Delta(q_1)$ to the origin. By the computation done in the proof of Proposition~\ref{periodcollapse}, we have~$\beta^{(p_1)}=(\frac{q_1}{3p_1},\frac{1}{3})$; hence~$\Delta=T_{(\frac{-q_1}{3p_1},\frac{-1}{3})}(\Delta(q_1))$.
\end{remark}
\paragraph{The rational case.} The following result shows that the Ehrhart theory of Markov triangles \nolinebreak{in standard~$p_1$-barycentric position} is particularly simple.
\begin{prop}\label{barycentreposition}
	All Markov triangles in standard~$p_1$-barycentric position are Ehrhart equivalent, having Ehrhart quasipolynomial with period~$3$.
\end{prop}

We note that the denominator of a Markov triangle representing~$(p_1,p_2,p_3)$ in standard~$p_1$-position is~$p_1p_2p_3$. The following result says that in this case there is strong, though not complete, period collapse.

\begin{prop}\label{periodcollapse}
	The Ehrhart quasipolynomial of a~$(p_1,p_2,p_3)$-triangle in standard~$p_1$-position has period a divisor of~$p_1$.
\end{prop}

The complete period collapse of the Fibonacci triangles proved in~\cite{gardiner_kleinman} is therefore a direct consequence of Proposition~\ref{periodcollapse}. Another direct implication is to the Markov triangles representing the triples lying on the~$2$-branch, which we call the \textit{Pell~triangles}.
\begin{corollary}\label{barnous}
	The Ehrhart quasipolynomial of a Pell triangle in standard~$2$-position has period~$2$.
\end{corollary}

\paragraph{The limiting case.} Take some Markov number~$a$ as the largest element of a Markov triple~$(b,a,c)$, and consider the bivalent subtree~$\Lambda(b,a,c)$. We study the evolution of the representative Markov triangles while going down the~$B_c$-branch of~$\Lambda(b,a,c)$: fix a companion number~$q\equiv 3cb^{-1}$~mod~$a$ of~$a$ with respect to~$(a,b,c)$, and consider the triple~$(a,b_n,c_n)$ obtained by first mutating~$(a,b,c)$ at~$b$ and then applying~$n-1$ increasing mutations preserving~$a$, with~$b_0\defeq b$,~$c_0\defeq c$, and~$b_{n+1}\defeq c_n$, $c_{n+1}\defeq3ac_n-b_n$. Let~$\Delta^a_n\defeq\Delta_n^a(q,c)$ be the Markov triangle representing~$(a,b_n,c_n)$ in standard~$a$-position with the~$b_n$-vertex on the~$x$-axis, and let~$\Delta^{a,\beta}_n\defeq \Delta^{a,\beta}_n(q,c)$ be the Markov triangle representing~$(a,b_n,c_n)$ in standard~$a$-barycentric position obtained by a translation of~$\Delta_n^a$. 

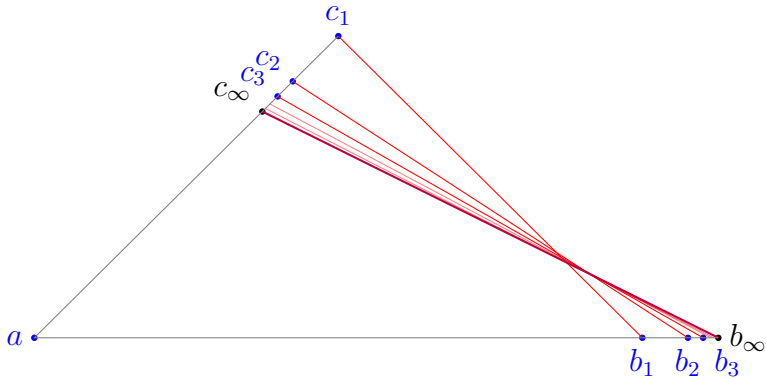
\begin{figure}[h]
	\centering
	\begin{tikzpicture}[scale=2]
		\filldraw [blue] (0,0) circle [radius=0.5pt] node[anchor=east] {$a$};
		
		\filldraw [blue] (4,0) circle [radius=0.5pt] node[anchor=north] {$b_1$};
		\filldraw [blue] (2,2) circle [radius=0.5pt] node[anchor=south]{$c_1$};
		\draw[color=red](4,0)--(2,2);
		
		\filldraw [blue] (4.3,0) circle [radius=0.5pt] node[anchor=north] {$b_2$};
		\filldraw [blue] (1.7,1.7) circle [radius=0.5pt] node[anchor=south east]{$c_2$};
		\draw[color=red](4.3,0)--(1.7,1.7);
		
		\filldraw [blue] (4.4,0) circle [radius=0.5pt] node[anchor=north west] {$b_3$};
		\filldraw [blue] (1.6,1.6) circle [radius=0.5pt] node[anchor=south east]{$c_3$};
		\draw[color=red](4.4,0)--(1.6,1.6);
		
		\draw[color=red,opacity=0.5](4.45,0)--(1.55,1.55);
		
		\draw[color=red, opacity=0.5](4.47,0)--(1.52,1.52);
		
		\filldraw [black] (4.5,0) circle [radius=0.5pt] node[anchor=west] {$b_\infty$};
		\filldraw [black] (1.5,1.5) circle [radius=0.5pt] node[anchor=south east]{$c_\infty$};
		\draw[color=purple, thick](4.5,0)--(1.5,1.5);
		
		\draw[color=gray] (4.5,0) -- (0,0) -- (2,2);
	\end{tikzpicture}
	\caption{Exaggerated representation of a sequence of Markov triangles in standard~$a$-position}
	\label{limitingtriangle}
\end{figure}

\begin{remarkqed}
	Going down the~$B_b$-branch of the bivalent tree instead of the~$B_c$-~branch yields a second pair of sequences of representative Markov triangles~$\Delta_n^a(-q,b)$ and~$\Delta_n^{a,\beta}(-q,b)$, both of which are constructed simply by swapping the letters~$b$ and~$c$ in the definitions of their~$B_c$-branch homologues. For the sake of clarity we will restrict our attention in the rest of the paper to the sequences~$\Delta_n^a(q,c)$ and~$\Delta_n^{a,\beta}(q,c)$, but note that the same results hold for the sequences~$\Delta_n^{a,\beta}(-q,b)$ and~$\Delta_n^{a,\beta}(-q,b)$.
\end{remarkqed}
Let us first investigate the sequence~$(\Delta_n^a)$ of triangles in standard~$a$-position; we note that~$\Delta_{n+1}^a$ is obtained by geometrically mutating~$\Delta_n^a$ at the~$b_n$-vertex for all~$n\geq 2$. The primitive vectors~$v_a^{(1)}\defeq v_{a,b_n}=(1,0)$ and~$v_a^{(2)}\defeq v_{a,c_n}=(aq-1,a^2)$ are preserved by the geometric mutations involved in the definition of~$\Delta_n^a$, and hence~$\Delta_n^a$ is entirely defined by the affine lengths of the~$[a,b_n]$- and~$[a,c_n]$-edges, namely~$\frac{c_n}{ab_n}$ and~$\frac{b_n}{ac_n}$ respectively, as illustrated in Figure~\ref{limitingtriangle}. We note that the ratio~$\frac{ab_n}{c_n}$ has a nice limiting behaviour,
\[\lim_{n\rightarrow\infty}\frac{ab_n}{c_n}=\frac{1}{\lambda(a)}\]
where~$\lambda(a)\defeq\frac{3+\sqrt{9-4/a^2}}{2}$ is the Lagrange number corresponding to~$a$, see~\nolinebreak{\cite[Section 3]{brendelschlenk}}\footnote{Brendel and Schlenk define~$\lambda(a)$ as~$\sqrt{9-4/a^2}$. We use an alternative definition in the present paper for the sake of clarity.}. We may then compute the limits of the affine lengths of the~$[a,b_n]$- and~$[a,c_n]$-edges
\begin{equation}\label{limitlength}
\lim_{n\rightarrow\infty}\frac{c_n}{ab_n}=\lambda(a)\qquad\text{and}\qquad \lim_{n\rightarrow\infty}\frac{b_n}{ac_n}=\frac{1}{a^2\lambda(a)}
\end{equation}
which allows us to define the triangle~$\Delta_\infty^a$ with~$a$-vertex at the origin,~$b_\infty$-vertex at an affine distance~$\lambda(a)$ from~$a$ along~$v_a^{(1)}$, and~$c_\infty$-vertex at an affine distance~$\frac{1}{a^2\lambda(a)}$ from~$a$ along~$v_a^{(2)}$.~$\Delta_\infty^a$ is in fact the limit of the Markov triangle sequence in question.

\begin{prop}\label{hausdorffconv}
	$(\Delta_n^a)$ Hausdorff converges to~$\Delta_\infty^a$.
\end{prop}

By Proposition~\ref{periodcollapse}, the Ehrhart function of~$\Delta_n^a$ is a quasipolynomial with period a divisor of~$a$, for all~$n\in\N$. It turns out that this period collapse also occurs in the limiting case:

\begin{theorem}\label{limittriangle}
	The~$a$-dilation of the limiting irrational triangle~$\Delta_\infty^a$ is Ehrhart equivalent to the~$a$-dilation of~$\Delta_n^a$ for any~$n\in\N$. In particular,~$\Delta_\infty^a$ is pseudo-rational with period~$a$.
\end{theorem}

This period collapse is not a priori a result of Hausdorff convergence, as pointed out in Remark~\ref{hausdorffnotehrhart}. Let us now consider the sequence~$(\Delta_n^{a,\beta})$ of Markov triangles in standard~$a$-barycentric position. It follows clearly from Proposition~\ref{hausdorffconv} and the fact that~$\Delta_n^{a,\beta}$ is the~$(\frac{-a}{3q},\frac{-1}{3})$-translation of~$\Delta_n^a$ by Remark~\ref{translatebarycentre}, that~$(\Delta_n^{a,\beta})$ Hausdorff converges to the triangle~$\Delta_\infty^{a,\beta}$ defined as the~$(\frac{-a}{3q},\frac{-1}{3})$-translation  of~$\Delta_\infty^a$.

\begin{theorem}\label{limitingbarycentre}
	The limiting irrational triangle~$\Delta_\infty^{a,\beta}$ in standard~$a$-barycentric position is Ehrhart equivalent to~$\Delta_n^{a,\beta}$ for any~$n\in\N$. In particular,~$3\Delta_\infty^{a,\beta}$ is pseudo-integral.
\end{theorem}


\paragraph{Open problem.} Recall from Example~\ref{fibonacci} the triangle~$\mathcal{F}_{(b/c,c/b)}$ for integers~$b,c\in\N$ with one vertex at the origin and the two other vertices at~$\frac{b}{c}(0,1)$ resp.\ $\frac{c}{b}(1,0)$. Cristofaro-Gardiner and Kleinman show in \cite[Theorem 1.6]{gardiner_kleinman} that~$\mathcal{F}_{(b/c,c/b)}$ is pseudo-integral if and only if~$b$ and~$c$ are consecutive odd-indexed Fibonacci numbers. This motivates the following problem: take coprime integers~$a,q\in\N$. For integers~$b,c\in\N$, define the triangle~$\mathcal{F}^{(p,q)}_{(b/c,c/b)}$ with one vertex at the origin and the two other vertices at~$\frac{b}{c}(aq-1,a^2)$ resp.\ $\frac{c}{b}(1,0)$. Proposition~\ref{periodcollapse} shows that~$\mathcal{F}^{(p,q)}_{(b/c,c/b)}$ is pseudo-integral if~$(a,b,c)$ is a Markov triple and~$q\equiv 3cb^{-1}$~mod~$a$, since in this case~$\mathcal{F}^{(p,q)}_{(b/c,c/b)}$ is the~$a$-dilate of a Markov triangle representing~$(a,b,c)$ in standard~$a$-position. This generalizes one direction of~\cite[Theorem 1.6]{gardiner_kleinman}. Is the converse true, i.e.\ do there exist pseudo-integral triangles~$\mathcal{F}^{(p,q)}_{(b/c,c/b)}$ which are not Markov? A similar question may be asked in the barycentric setting: Proposition~\ref{barycentreposition} says that the~$\frac{3}{a}$-dilate of the~$(\frac{-q}{3a},\frac{-1}{3})$-translate of the triangle~$\mathcal{F}^{(p,q)}_{(b/c,c/b)}$ is pseudo-integral if~$(a,b,c)$ is a Markov triple and~$q\equiv 3cb^{-1}$~mod~$a$. Is the converse true?
	\subsection{Connections to symplectic geometry}\label{section_5}
	The rational triangles whose Ehrhart theory we study in this paper are special in two ways: they admit complete or at least strong period collapse, and their irrational limit triangles also have this property. This Ehrhart-peculiarity may be a shadow of their relevance in symplectic geometry.

Indeed, our rational triangles arise as bases of almost toric fibrations of the complex projective plane~$\mathbb{CP}^2$, see~\cite{vianna},~\cite{evans2023lectures}. Furthermore, these triangles are the moment map images of exactly those weighted projective planes that admit smoothings to~$\mathbb{CP}^2$, see~\cite{Hacking_Prokhorov_2010}. On the other hand, these triangles also describe certain rational homology ellipsoids, that in the special case of Fibonacci triangles are usual ellipsoids. As a consequence, for every Markov number, the problem of symplectically embedding rational homology ellipsoids into a minimal~$\mathbb{CP}^2$ leads to an infinite staircase, see~\cite{adaloglou2025markovstaircases}. 

For Fibonacci triangles, this relation between Ehrhart theory and symplectic embedding problems was noticed and studied  earlier in~\cite{gardiner_kleinman} and~\cite{cristofarogardiner2015newexamplesperiodcollapse}, and the corresponding problem of symplectically embedding usual ellipsoids into~$\mathbb{CP}^2$, but also into the 4-ball, led to the Fibonacci staircase in~\cite{mcduffschlenk}.
	\newpage
	\section{Proofs}\label{proofs}
	\subsection{Integral barycentre of a rational triangle}\label{section_21}
Let~$p,p_0\in\R^2$ be rational vectors,~$v\in\Z^2$ a primitive integral vector. We define the \textit{affine distance} from~$p$ to the affine line~$\ell\defeq\ell(p_0,v)$ through~$p_0$ spanned by~$v$ to be the absolute value of the determinant of~$p-p_0$ with~$v$, i.e.~$|\det(p-p_0,v)|$.

\begin{remarkqed}\label{distancechoice}
	The affine distance from~$p$ to~$\ell$ does not depend on the choice of~$p_0$: given another point~$p_0'\in\ell$, we have~$$\det(p-p_0',v)=\det(p-p_0+p_0-p_0',v)=\det(p-p_0,v)$$ since~$p_0-p_0'$ is a scalar multiple of~$v$.
\end{remarkqed}
In classical geometry, the bisector of an angle~$\alpha$ in the plane is defined as the set of all points lying at an equal distance from both lines spanning~$\alpha$. The following result asserts that a similar definition holds true for integral bisectors.
\begin{lemma}\label{bisectordistance}
	Let~$u,v\in\Z^2$ be linearly independent primitive integral vectors. A~vector~$w\in\R^2$ is a scalar multiple of~$u+v$ if and only if~$\det(w,u)=\det(v,w)$.
\end{lemma}
\begin{proof}
	Assume first that~$w=\lambda(u+v)$ for some~$\lambda\in\R$. Then
	\[\det(w,u)=\det(\lambda(u+v),u)=\lambda\det(v,u)=\det(v,\lambda(u+v))=\det(v,w).\]
	Conversely assume that~$w$ is not a multiple of~$u+v$. Since~$u$ and~$v$ span~$\R^2$, then we may write~$w=\lambda_1u+\lambda_2v$ for some scalars~$\lambda_1\neq\lambda_2$. Then 
	\[\det(w,u)=
		\lambda_2\det(v,u)
		\neq\lambda_1\det(v,u)
		=\det(v,w).\qedhere\]
\end{proof}

Let~$\Delta\subset\R^2$ be a rational triangle with vertices~$p_1,p_2,p_3$. Unwinding the definitions, Lemma~\ref{bisectordistance} says that the integral bisector of the angle at~$p_1$ is the set of points at equal affine distance from the~$[p_1,p_2]$- and the~$[p_1,p_3]$-edges.
\begin{lemma}\label{barycentre}
	The three integral bisectors of~$\Delta$ intersect at a single point.
\end{lemma}
\begin{proof}
	Let~$\ell_i$ be the integral bisector of the angle at~$p_i$ (for~$i\in\Z_3$), and let~$\beta$ be the intersection of~$\ell_1$ and~$\ell_2$. Then we have
	\begin{align*}
		|\det(\beta-p_3,v_{1,3})|&=|\det(\beta-p_1,v_{1,3})|=|\det(v_{1,2},\beta-p_1)|\\
		&=|\det(v_{1,2},\beta-p_2)|=|\det(\beta-p_2,v_{2,3})|=|\det(\beta-p_3,v_{2,3})|
	\end{align*}
	where the first, third and fifth equalities follow from Remark~\ref{distancechoice}, while the second and fourth equalities follow from the definition of~$\beta$ and Lemma~\ref{bisectordistance}. We deduce, also by Lemma~\ref{bisectordistance}, that~$\beta-p_3$ is a scalar multiple of~$v_{3,1}+v_{3,2}$, i.e.\ it lies on~$\ell_3$.
	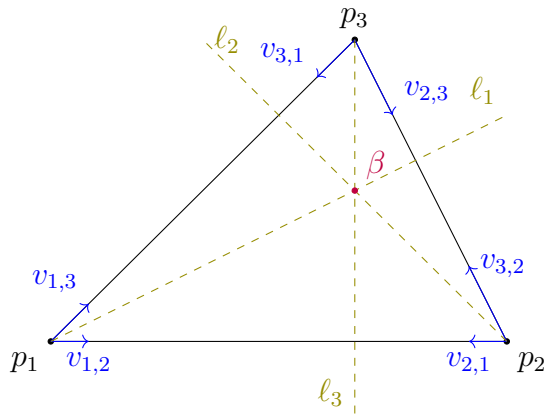
\begin{figure}[h]
	\centering
	\begin{tikzpicture}[scale=2]
		\filldraw [black] (0,0) circle [radius=0.5pt] node[anchor=north east] {$p_1$};
		\filldraw [black] (3,0) circle [radius=0.5pt] node[anchor=north west] {$p_2$};
		\filldraw [black] (2,2) circle [radius=0.5pt] node[anchor=south]{$p_3$};
		\draw[color=black](0,0)--(3,0)--(2,2)--cycle;
		
		\draw[color=blue,->](0,0)--(0.25,0.25) node[anchor=south east] {$v_{1,3}$};
		\draw[color=blue,->](2,2)--(1.75,1.75) node[anchor=south east] {$v_{3,1}$};
		\draw[color=blue,->](0,0)--(0.25,0) node[anchor=north] {$v_{1,2}$};
		\draw[color=blue,->](3,0)--(2.75,0) node[anchor=north] {$v_{2,1}$};
		\draw[color=blue,->](3,0)--(2.75,0.5) node[anchor=west] {$v_{3,2}$};
		\draw[color=blue,->](2,2)--(2.25,1.5) node[anchor=south west] {$v_{2,3}$};
		
		\draw[color=olive,dashed] (0,0)--(3,1.5) node[anchor=south east] {$\ell_1$};
		\draw[color=olive,dashed] (2,2)--(2,-0.5) node[anchor=south east] {$\ell_3$};
		\draw[color=olive,dashed] (3,0)--(1,2) node[anchor=west] {$\ell_2$};
		
		\filldraw [purple] (2,1) circle [radius=0.5pt] node[anchor=south west]{$\beta$};
	\end{tikzpicture}
	\caption{Proof of Lemma~\ref{barycentre}}
	\label{barycentrepic}
	\end{figure}\qedhere
\end{proof}
As mentioned in~\S\ref{section_4}, we call the intersection of the integral bisectors of~$\Delta$, the \textit{integral barycentre} of~$\Delta$.

\begin{remark}
	The integral barycentre of~$\Delta$ is equivalently defined as the point lying at an equal affine distance from the three edges of~$\Delta$.
\end{remark}

\subsection{On~model Markov triangles}\label{section_22}
We introduce an alternative version of Markov triangles, which turns out to be equivalent to Definition~\ref{defmarkov}.
\begin{definitionqed}\label{pseudo}
	Let~$(p_1,p_2,p_3)$ be a Markov triple. We define a \textit{pseudo-Markov triangle}~$\Delta_{(p_1,p_2,p_3)}$ representing~$(p_1,p_2,p_3)$ to be any triangle (with vertices denoted by~$p_1,p_2,p_3$) satisfying the following conditions
	\begin{itemize}
		\item the edge opposite to~$p_i$ has affine length~$\frac{p_i}{p_{i+1}p_{i+2}}$ (where~$i\in\Z_3$),
		\item the determinant of the angle at the~$p_i$-vertex is~$p_i^2$.\qedhere
	\end{itemize}
\end{definitionqed}

To prove the equivalence of Definitions~\ref{defmarkov} and~\ref{pseudo}, we need the following two facts about companion numbers of Markov numbers.
\begin{lemma}\label{mutationcompanion}
	Consider a Markov triple~$(p_1,p_2,p_3)$ and take a companion number~$q_1$ of~$p_1$ with respect to~$(p_1,p_2,p_3)$. Then~$q_1$~and~$p_1$ are coprime. Moreover,~$q_1$ is also a companion of~$p_1$ with respect to~$(p_1,\widehat{p_2},p_3)$.
\end{lemma}
\begin{proof}
	First we show that~$p_1$ is relatively prime to~$3$: assume by way of contra\-diction that~$3\mid p_1$, so that the Markov equation implies~$p_2^2+p_3^2\equiv 0$~mod~$3$. Since~$x^2\equiv 0$ or~$x^2\equiv 1$~mod~$3$ for all~$x$, it follows that~$p_2^2\equiv p_3^2\equiv 0$~mod~$3$ and hence~$3$ divides all three elements of the triple~$(p_1,p_2,p_3)$, contra\-dicting their pairwise coprimality.\linebreak Thus~$p_1$ is coprime to any number congruent to~$\pm 3p_3p_2^{-1}$~mod~$p_1$, in particular $q_1$. In addition, we have
		\begin{align*}
			q_1\equiv \pm 3p_2p_3^{-1}&\equiv \mp 3(-p_2)p_3^{-1}&\text{mod }p_1\\
			&\equiv \mp 3(3p_1p_3-p_2)p_3^{-1}&\text{mod }p_1\\
			&\equiv \mp 3\widehat{p_2}p_3^{-1}&\text{mod }p_1
		\end{align*}
	which means that $q_1$ is a companion of~$p_1$ with respect to~$(p_1,\widehat{p_2},p_3)$.
\end{proof}

The following proposition gives three ways to think about the rational triangles studied in this paper.
\begin{prop}\label{pseudomarkov}
	Let~$\Delta\subset\R^2$ be a triangle and consider a Markov triple~$(p_1,p_2,p_3)$. The following statements are equivalent.
	\begin{enumerate}
		\item~$\Delta$ is a Markov triangle representing~$(p_1,p_2,p_3)$.
		\item~$\Delta$ is a pseudo-Markov triangle representing~$(p_1,p_2,p_3)$.
		\item~$\Delta$ is~$\R^2\rtimes\operatorname{GL}(2,\Z)$-congruent to a triangle having vertices
		\begin{equation}\label{bijour}
			\Bigg\{(0,0),\Big(\frac{(p_1q_1-1)p_2}{p_1p_3},\frac{p_1p_2}{p_3}\Big),\Big(\frac{p_3}{p_1p_2},0\Big)\Bigg\}
		\end{equation}
		for some companion number~$q_1\equiv 3p_3p_2^{-1}$~mod~$p_1$ of~$p_1$.
	\end{enumerate}
\end{prop}
\begin{proof}
	Let~$v_{i,j}$ be the primitive vector along the~$[p_i,p_j]$-edge of~$\Delta$. Fix a companion number~$q_1\equiv p_3p_2^{-1}$~mod~$p_1$ of~$p_1$, and let~$\Delta_0$ be the triangle with~$p_1$-,~$p_3$-, and~$p_2$-vertices given by \eqref{bijour}.
	
	(2)$\Rightarrow$(3): Assume~$\Delta$ is a pseudo-Markov triangle representing~$(p_1,p_2,p_3)$. By~Remark~\ref{congruentmarkov}, we may assume that the~$p_1$-vertex of~$\Delta$ is at the origin and \mbox{the~$[p_1,p_2]$-edge} lies on the~$x$-axis; this edge is then the interval~$[0,\frac{p_3}{p_1p_2}]$. Hence we already have the coordinates of the~$p_1$- and~$p_2$-vertices, and turn to those of the~$p_3$-vertex. We have~$v_{1,2}=(1,0)$. By the determinant property of~$\Delta$ at the~$p_1$-vertex,\linebreak we have~$|y_{v_{1,3}}|=|\det(v_{1,3},(1,0))|=p_1^2$, so that~$y_{v_{1,3}}=\pm p_1^2$. Applying the~$\operatorname{GL}(2,\Z)$-~transformation~$(x,y)\mapsto(x,-y)$ to~$\Delta$ if necessary, we can assume that~$y_{v_{1,3}}=p_1^2$. Arguing similarly at the~$p_2$-vertex yields~$y_{v_{2,3}}=p_2^2$, so~$v_{1,3}=(a_1,p_1^2)$ and~$v_{2,3}=(a_2,p_2^2)$ for some integers~$a_1,a_2\in\Z$. The determinant at the~$p_3$-vertex is
	\begin{equation}\label{bornos}
		p_3^2=|a_1p_2^2-a_2p_1^2|=a_1p_2^2-a_2p_1^2.
	\end{equation}
	where we used the fact that~$a_1p_2^2-a_2p_1^2=\det(v_{3,1},v_{3,2})>0$ because of our choice of orientation for~$\Delta$: counterclockwise through~$p_1$, then~$p_2$, then~$p_3$. From this and the Markov equation, we have
	\[a_1p_2^2=p_3^2+a_2p_1^2=3p_1p_2p_3-p_1^2-p_2^2+a_2p_1^2\equiv -p_2^2 \qquad \text{mod }p_1\]
	and hence~$a_1\equiv -1$~mod~$p_1$ since~$p_1$ and~$p_2$ are coprime. In other terms,~$a_1+1$ is divi\-sible by~$p_1$, allowing us to define the integer~$q_1'\defeq\frac{a_1+1}{p_1}$, for which~$v_{1,3}=(p_1q_1'-1,p_1^2)$. We have
	\begin{equation*}
		q_1'p_2^2=\frac{a_1p_2^2+p_2^2}{p_1}=\frac{3p_1p_2p_3-p_1^2+a_2p_1^2}{p_1}=3p_2p_3+p_1(a_2-1)
	\end{equation*}
	so that~$q_1'p_2^2\equiv 3p_2p_3$~mod~$p_1$. Multiplying through by~$(p_2^{-1})^2$ yields~$q_1'\equiv 3p_3p_2^{-1}$~mod~$p_1$. The coordinates of the~$p_3$-vertex are then given by
	\[\Big(\frac{(p_1q_1'-1)p_2}{p_1p_3},\frac{p_1p_2}{p_3}\Big)\]
	because the~$p_3$-vertex belongs to the line~$y=\frac{p_1^2}{p_1q_1'-1}x$ spanned by~$v_{1,3}$, and is at an affine distance~$\frac{p_2}{p_1p_3}$ from the origin. Since~$q_1\equiv q_1'$ mod~$p_1$, we are done by Remark~\ref{equimarkov}.
	
	(3)$\Rightarrow$(2): Consider the triangle~$\Delta_0$. We have
	\[q_1p_2-3p_3\equiv 3p_3-3p_3\equiv 0\qquad\text{mod }p_1\]
	so that we can define the integer~$q_2\defeq\frac{q_1p_2-3p_3}{p_1}$, which is in particular a companion number of~$p_2$, namely~$q_2\equiv- 3p_3p_1^{-1}$~mod~$p_2$.
	We first verify that~$\Delta_0$ satisfies the edge property for being pseudo-Markov.
	\begin{itemize}
		\item Dividing the vector from the~$p_2$-vertex to the~$p_3$-vertex by~$\frac{p_1}{p_2p_3}$ yields~$(p_2q_1+1,p_2^2)$, which is primitive.
		\item The vector~$\frac{p_1p_3}{p_2}\big(\frac{(p_1q_1-1)p_2}{p_1p_3},\frac{p_1p_2}{p_3}\big)=(p_1q_1-1,p_1^2)$ is primitive.
		\item The affine length of the~$[p_1,p_2]$-edge is of course~$\frac{p_3}{p_1p_2}$.
	\end{itemize}
	We then verify the angle property.
	\begin{itemize}
		\item At~$p_1$:~$v_{1,2}=(1,0)$ and~$v_{1,3}=(p_1q_1-1,p_1^2)$, so~$|\det(v_{1,2},v_{1,3})|=p_1^2$.
		
		\item At~$p_2$:~$v_{2,1}=(-1,0)$~and~$v_{2,3}=(p_2q_2+1,p_2^2)$, so~$|\det(v_{2,1},v_{2,3})|=p_2^2$.
		
		\item At~$p_3$:~$v_{3,1}=(1-p_1q_1,-p_1^2)$ ~and~$v_{3,2}=(-p_2q_2-1,-p_2^2)$, so
		\begin{align*}
			|\det(v_{3,1},v_{3,2})|&=|-p_2^2(1-p_1q_1)+p_1^2(-p_2q_2-1)|\\
			&=|-p_1^2-p_2^2+p_2^2p_1q_1-p_1^2p_2q_2|\\
			&=|-p_1^2-p_2^2+p_1p_2(p_2q_1-p_1q_2)|\\
			&=|-p_1^2-p_2^2+3p_1p_2p_3|=|p_3^2|=p_3^2.
		\end{align*}
	\end{itemize}
	
	(3)$\Leftrightarrow$(1): In the following we assume without loss of generality that~$p_2\geq \max(p_1,p_3)$, so that~$(p_1,p_3,\widehat{p_2})$ is the parent of~$(p_1,p_2,p_3)$ in the Markov tree; in particular~$(p_1,p_3,\widehat{p_2})$ lives in the~$(n-1)$-st generation of the tree. We prove the equivalence (3)$\Leftrightarrow$(1) by induction on the generation~$n$ of the triple~$(p_1,p_2,p_3)$. For~$n=1$, the triple is~$(1,1,1)$ and in this case the representing Markov triangles are \textit{defined} as the pseudo-Markov triangles, so we have (1)$\Leftrightarrow$(2). Since (2)$\Leftrightarrow$(3), we are done.\\
	Assume that (3) and (1) are equivalent for all~$(n-1)$-st generation triples of the Markov tree. Consider the triangle~$\widehat{\Delta}$ with~$p_1$-,~$\widehat{p_2}$-, and~$p_3$-vertices
	\begin{equation*}
	\Bigg\{(0,0),\Big(\frac{(p_1q_1-1)p_3}{p_1\widehat{p_2}},\frac{p_1p_3}{\widehat{p_2}}\Big),\Big(\frac{\widehat{p_2}}{p_1p_3},0\Big)\Bigg\}.
	\end{equation*}
	As shown in the proof of Lemma~\ref{mutationcompanion}, we have~$q_1\equiv -3p_3\widehat{p_2}^{-1}\equiv 3\widehat{p_2}p_3^{-1}$~mod~$p_1$, so~$\widehat{\Delta}$ satisfies property 3; by the inductive hypothesis, the Markov triangles re\-presenting~$(p_1,p_3,\widehat{p_2})$ are precisely the elements of the~$\R^2\rtimes\operatorname{GL}(2,\Z)$-orbit of~$\widehat{\Delta}$.
	
	\begin{claim}
		$\widehat{\Delta}$ is the image of~$\Delta_0$ under the geometric mutation~$M_{p_2}$ at \nolinebreak{the~$p_2$-vertex}.	
	\end{claim}
	
	The equivalence (3)$\Leftrightarrow$(1) follows from this claim:~$\Delta$ satisfies property~3 if and only if there exists a map~$\varphi\in\R^2\rtimes\operatorname{GL}(2,\Z)$ for which
	\[\varphi(\Delta)=\Delta_0\qquad\Leftrightarrow\qquad M_{p_2}(\varphi(\Delta))=M_{p_2}(\Delta_0)=\widehat{\Delta}\qquad\Leftrightarrow\qquad\varphi(\Delta)=M_{\widehat{p_2}}(\widehat{\Delta})\]
	i.e. $\varphi(\Delta)$ is the geometric mutation of a Markov triangle representing~$(p_1,p_3,\widehat{p_2})$~at~$\widehat{p_2}$ and is therefore a Markov triangle representing~$(p_1,p_2,p_3)$ itself. By Remark~\ref{congruentmarkov}, this means that~$\Delta$ is a Markov triangle representing~$(p_1,p_2,p_3)$, i.e. $\Delta$ satisfies property~1.
	
	\noindent\textit{Proof of the claim.} Upon setting~$q_2\defeq \frac{q_1p_2-3p_3}{p_1}$, the primitive vector~$v_{2,3}$ along the~$[p_2,p_3]$-edge of~$\Delta_0$ is given by~$(p_2q_2+1,p_2^2)$ as in the proof of~(3)$\Rightarrow$(2), and hence the primitive vector~$v_{p_2}$ along the integral bisector of the~$p_2$-vertex~$\ell_2$ is~$v_{p_2}=(q_2,p_2)$. We first verify that the~$\widehat{p_2}$-vertex in~$\widehat{\Delta}$ is the intersection of~$\ell_2$ with the~$[p_1,p_3]$-edge, i.e.\ that the~$\widehat{p_2}$-vertex satisfies the equations
	\[y=\frac{p_2}{q_2}\big(x-\frac{p_3}{p_1p_2}\big)\qquad\text{and}\qquad y=\frac{p_1^2}{p_1q_1-1}x.\]
	Solving this system for~$x$ yields
	\[x=\frac{p_3(p_1q_1-1)}{p_1\big(p_1(p_2q_1-p_1q_2)-p_2\big)}.\]
	Using the expressions for~$\widehat{p_2}$ and~$q_2$ in terms of~$p_1,p_2,p_3,q_1$, we get
	\[x=\frac{p_3(p_1q_1-1)}{p_1(3p_3p_1-p_2)}=\frac{p_3(p_1q_1-1)}{p_1\widehat{p_2}}\qquad\text{and}\qquad y=\frac{p_1^2}{p_1q_1-1}x=\frac{p_1p_3}{\widehat{p_2}}\] as claimed. We check that~$M_{p_2}$ takes the~$p_3$-vertex~of~$\Delta$ to the~$p_3$-vertex~of~$\widehat{\Delta}$. The shear~$\varphi_{v_{p_2}}$ is given by
	\begin{equation}\label{abouali}
		\varphi_{v_{p_2}}=\begin{pmatrix}
			1-x_{v_{p_2}}y_{v_{p_2}}&x_{v_{p_2}}^2\\
			-y_{v_{p_2}}^2&1+x_{v_{p_2}}y_{v_{p_2}}
		\end{pmatrix}=\begin{pmatrix}
			1-q_2p_2&q_2^2\\
			-p_2^2&1+q_2p_2
		\end{pmatrix}.
	\end{equation}
	We compute
	\[M_{p_2}\begin{pmatrix}
		\frac{(p_1q_1-1)p_2}{p_1p_3}\\[6pt]
		\frac{p_1p_2}{p_3}
	\end{pmatrix}=\begin{pmatrix}
	\frac{q_2\big(p_1p_2(q_2p_1-p_2q_1)+p_2^2\big)+p_2p_1q_1+q_2p_3^2-p_2}{p_1p_3}\\[6pt]
	\frac{p_1p_2^3(q_2p_1-q_1p_2)+p_2^2(p_1^2+p_2^2+p_3^2)}{p_1p_2p_3}
	\end{pmatrix}.\]
	Using the expressions for~$\widehat{p_2}$ and~$q_2$ again as well as the Markov equation, we simplify
	\begin{align*}
		q_2\big(p_1p_2(q_2p_1-p_2q_1)+p_2^2\big)&+p_2p_1q_1+q_2p_3^2-p_2\\
		&=q_2(-3p_1p_2p_3+p_2^2)+p_2p_1q_1+q_2p_3^2-p_2\\
		&=q_2(-p_1^2-p_3^2)+p_2p_1q_1+q_2p_3^2-p_2\\
		&=p_1(q_1p_2-q_2p_1)-p_2=3p_1p_3-p_2=\widehat{p_2}
	\end{align*}
	and
	\begin{align*}
		p_1p_2(q_2p_1-q_1p_2)+p_1^2+p_2^2+p_3^2=-3p_1p_2p_3+p_1^2+p_2^2+p_3^2=0.
	\end{align*}
	Hence we have
	\[M_{p_2}\begin{pmatrix}
		\frac{(p_1q_1-1)p_2}{p_1p_3}\\[6pt]
		\frac{p_1p_2}{p_3}
	\end{pmatrix}=\begin{pmatrix}
	\frac{\widehat{p_2}}{p_1p_3}\\[6pt]
	0
	\end{pmatrix}\]
	as desired.
\end{proof}
We observe that the geometric mutations of a Markov triangle preserve its integral barycentre. More precisely,
\begin{lemma}\label{mutatebarycentre}
	Let~$\Delta$ be a Markov triangle representing~$(p_1,p_2,p_3)$ and~$\beta$ the barycentre of~$\Delta$. Then~$\beta$ is also the barycentre of~$M_{p_i,\Delta}(\Delta)$ for~$i\in\Z_3$.
\end{lemma}
\begin{proof}
	Let~$\ell_i$ be the integral bisector at the~$p_i$-vertex of~$\Delta$.
	Without loss of generality, we fix the orientation of~$\Delta$ to be counter-clockwise through~$p_1,p_2,p_3$, and restrict our attention to the geometric mutation~$M_{p_2}\defeq M_{p_2,\Delta}$ at~$p_2$.
	Thus~$\widehat{\Delta}\defeq M_{p_2}(\Delta)$ is the triangle with vertices~$p_1$,~$M_{p_2}(p_3)$, and~$\widehat{p_2}\defeq\ell_2\cap[p_1,p_3]$. We denote the edges of~$\widehat{\Delta}$ by~$\widehat{e_{1,2}}\defeq[p_1,\widehat{p_2}]$,~$\widehat{e_{1,3}}\defeq[p_1,M_{p_2}(p_3)]$, and~$\widehat{e_{2,3}}\defeq[\widehat{p_2},M_{p_2}(p_3)]$ and let~$\widehat{v_{i,j}}$ be the primitive integral vector along~$\widehat{e_{i,j}}$. We observe that~$\widehat{v_{1,3}}=v_{1,2}$ and~$\widehat{v_{1,2}}=v_{1,3}$, so that the integral bisector~$\widehat{\ell_1}$ at the~$p_1$-vertex of~$\widehat{\Delta}$ is still~$\ell_1$. On the other hand, we claim that the integral bisector~$\widehat{\ell_2}$ at the~$\widehat{p_2}$-vertex is~$\ell_2$, because the points lying on~$\ell_2$ are at an equal affine distance from both edges~$\widehat{e_{1,2}}$ and~$\widehat{e_{2,3}}$. This property is trivially satisfied by the~$\widehat{p_2}$-vertex, so we only need to find one other point of~$\ell_2$ satisfying it. Consider the~$p_2$-vertex of~$\Delta$. The affine distance from~$p_2$ to~$\widehat{e_{2,3}}$ is given by the absolute value of the determinant of~$v\defeq M_{p_2}(p_3)-p_2$ with~$\widehat{v_{2,3}}$. Since~$v$ is the image of~$p_3-p_2$ under a shear, both vector have the same affine length, equal to~$\frac{p_1}{p_2p_3}$ by the edge property of~$\Delta$ as a pseudo-Markov triangle in Definition~\ref{pseudo}. Thus we have
	\[\operatorname{dist}(p_2,\widehat{e_{2,3}})=|\det(v,\widehat{v_{2,3}})|=\frac{p_1}{p_2p_3}|\det(v_{1,2},\widehat{v_{2,3}})|=\frac{p_2p_3}{p_2}\]
	where~$|\det(v_{2,1},\widehat{v_{2,3}})|=p_3^2$ by the angle property of~$\widehat{\Delta}$ as a pseudo-Markov triangle. Likewise the affine distance from~$p_2$ to~$\widehat{e_{1,2}}$ is given by
	\[\operatorname{dist}(p_2,\widehat{e_{1,2}})=|\det(v,\widehat{v_{1,2}})|=\frac{p_3}{p_1p_3}|\det(v_{1,2},\widehat{v_{1,3}})|=\frac{p_2p_3}{p_2}\]
	We conclude that~$\widehat{\ell_2}=\ell_2$. Therefore the integral barycentre~$\widehat{\beta}$ of~$\widehat{\Delta}$ is given by~$\widehat{\beta}=\widehat{\ell_1}\cap\widehat{\ell_2}=\ell_1\cap\ell_2=\beta$.
\end{proof}
	
	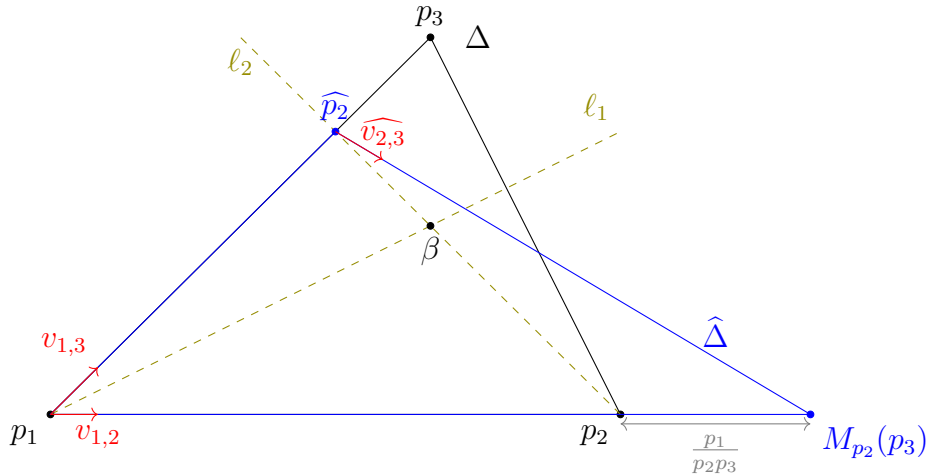
\begin{figure}[h]
		\centering
		\begin{tikzpicture}[scale=2.5]
			\filldraw [black] (0,0) circle [radius=0.5pt] node[anchor=north east] {$p_1$};
			\filldraw [black] (3,0) circle [radius=0.5pt] node[anchor=north east] {$p_2$};
			\filldraw [black] (2,2) circle [radius=0.5pt] node[anchor=south]{$p_3$};
			\draw[color=black](0,0)--(3,0)--(2,2)--cycle;
			
			\draw[color=olive,dashed] (0,0)--(3,1.5) node[anchor=south east] {$\ell_1$};
			\draw[color=olive,dashed] (3,0)--(1,2) node[anchor=north] {$\ell_2$};
			
			\filldraw [black] (2,1) circle [radius=0.5pt] node[anchor=north]{$\beta$};
			
			\draw [black] (2.25,2) node{$\Delta$};
			\draw [blue] (3.5,0.3) node[anchor=south]{$\widehat{\Delta}$};
			
			\filldraw [blue] (1.5,1.5) circle [radius=0.5pt] node[anchor=south] {$\widehat{p_2}$};
			\filldraw [blue] (4,0) circle [radius=0.5pt] node[anchor=north west] {$M_{p_2}(p_3)$};
			\draw[color=gray,<->] (3,-0.05)--(4,-0.05);
			\draw[color=gray](3.5,-0.05)node[anchor=north]{$\frac{p_1}{p_2p_3}$};
			\draw[color=blue](0,0)--(4,0)--(1.5,1.5)--cycle;
			
			
			\draw[color=red,->](1.5,1.5)--(1.75,1.35) node[anchor=south] {$\widehat{v_{2,3}}$};
			\draw[color=red,->](0,0)--(0.25,0.25) node[anchor=south east] {$v_{1,3}$};
			\draw[color=red,->](0,0)--(0.25,0) node[anchor=north] {$v_{1,2}$};
		\end{tikzpicture}
		\caption{Preservation of the integral barycentre under geometric mutation}
		\label{preservebary}
	\end{figure}

\begin{remark}\label{distancebarycentre}
	It follows directly from Lemma~\ref{mutatebarycentre} and its proof that the affine distance from~$\beta$ to any edge of~$\Delta$ is preserved by the geometric mutation of~$\Delta$ at any vertex.
\end{remark}

\subsection{On the Ehrhart theory of rational Markov triangles}\label{section_23}
\begin{proof}[Proof of Proposition~\ref{brunain}]
	Let~$P\subset\R^2$ be a polygon and~$v$ an integral vector. Consider the half-shear~$\varphi\defeq\varphi_v^{1/2}:\R^2\rightarrow\R^2$ with respect to~$v$, which preserves the integer lattice as mentioned in Remark~\ref{shearpreserveslattice}. Then for any point~$x\in P$ and integer dilation factor~$t\in\N\cup\{0\}$,~$tx\in\Z^2$ if and only if~$t\varphi(x)=\varphi(tx)\in\Z^2$, so that~$\varphi$ defines a bijection between the integer lattice points in the~$t$-dilate of~$P$ and in the~$t$-dilate of~$\varphi(P)$, for each~$t\in\N\cup\{0\}$. Since~$P$ has finitely many integer lattice points, we get 
	\[L_P(t)=\#(tP)\cap\Z^2=\#\big(t\varphi(P)\cap\Z^2\big)=L_{\varphi(P)}(t).\qedhere\]
\end{proof}

\begin{proof}[Proof of Proposition~\ref{barycentreposition}]
	Let~$(p_1,p_2,p_3)$ be a Markov triple and consider a re\-presentative Markov triangle~$\Delta^\beta$ in standard~$p_1$-barycentric position. For each ~$i=1,2,3$, the integral barycentre~$\beta$ of~$\Delta$ lies on the integral bisector~$\ell_i$ of the~$p_i$-vertex, so that the position vector~$p_i=p_i-\beta$ is parallel to~$\ell_i$. It follows that the half-shear~$\varphi_{v_i}^{1/2}$ with respect to the primitive vector~$v_i$ along~$\ell_i$ is invariant under translation along~$p_i$, and hence the geometric mutation~$M_{p_i}=T_{p_i}\circ\varphi_{v_i}^{1/2}\circ T_{-p_i}=\varphi_{v_i}^{1/2}$ is the half-shear itself, so that~$M_{p_i}(\Delta^\beta)$ and~$\Delta^\beta$ are Ehrhart equivalent by Proposition~\ref{brunain}.
\end{proof}

For the proof of Proposition~\ref{periodcollapse}, we need the following two results on the effect of geometric mutations on Markov triangles in standard~$p$-position.
\begin{lemma}\label{ehrhartmutate}
	Let~$(p_1,p_2,p_3)$ be a Markov triple and~$\Delta$ a representing triangle in standard~$p_1$-position.
	Let~$\widehat{\Delta}$ be the image of~$\Delta$ under its geometric mutation at any vertex. Then~$p_1\Delta$ and~$p_1\widehat{\Delta}$ are Ehrhart equivalent.
\end{lemma}
\begin{proof}
	In the following we assume that the~$p_2$-vertex lies on the~$x$-axis. The result is trivial for the geometric mutation of the~$t$-dilate of~$\Delta$ at~$tp_1$, since~$tp_1$ is the origin and hence~$M_{tp_1,t\Delta}=M_{(0,0),t\Delta}$ is the half-shear along~$v_1$, and in particular preserves the integer lattice. Thus~$M_{p_1,\Delta}(\Delta)$ and~$\Delta$ are Ehrhart equivalent.
	
	We next consider the geometric mutation of the~$t$-dilate of~$\Delta$ at~$tp_2$, defined as
	\[M_{tp_2, t\Delta}=T_{(\frac{tp_3}{p_1p_2},0)}\circ\varphi_{v_{p_2}}^{1/2}\circ T_{(\frac{-tp_3}{p_1p_2},0)}\]
	as described in the proof of Proposition~\ref{pseudomarkov}, where the half-shear~$\varphi_{v_{p_2}}^{1/2}$ is given on~$H_{v_{p_2}}^+$ in matrix form by \eqref{abouali}. Take an integer point~$(x,y)\in\Z^2$. If~$(x,y)$ lies to the left of the integral bisector~$\ell(tp_2,v_{p_2})$ at the~$tp_2$-vertex of~$t\Delta$, 
	then~$M_{tp_2,t\Delta}(x,y)=(x,y)$. Otherwise,
	
	\vspace{-15mm}
	
	\begin{align*}
		M_{tp_2,t\Delta}(x,y)&=\begin{pmatrix}
			\frac{tp_3}{p_1p_2}\\[6pt]
			0
		\end{pmatrix}+\begin{pmatrix}
			1-p_2q_2&q_2^2\\[6pt]
			-p_2^2&1+p_2q_2
		\end{pmatrix}\begin{pmatrix}
			x-\frac{tp_3}{p_1p_2}\\[6pt]
			y
		\end{pmatrix}\\
		&=\begin{pmatrix}
			x(1-p_2q_2)+y{q_2}^2+q_2p_3\frac{t}{p_1}\\[6pt]
			x(-p_2^2)+y(1+p_2q_2)+p_2p_3\frac{t}{p_1}
		\end{pmatrix}\\
		&=\begin{pmatrix}
			k_1x+k_1'y+k_1''\frac{t}{p_1}\\[6pt]
			k_2x+k_2'y+k_2''\frac{t}{p_1}
		\end{pmatrix}
	\end{align*}
	for some integers~$k_1,k_1',k_1'',k_2,k_2',k_2''\in\Z$. Hence~$M_{tp_2,t\Delta}(\Z^2)\subset\Z^2$ when~$p_1$ divides~$t$; the sublattice~$M_{tp_2,t\Delta}(\Z^2)$ of~$\Z^2$ has index~$\det(M_{tp_2,t\Delta})=1$, so that~$M_{tp_2,t\Delta}(\Z^2)=\Z^2$.
	Thus~$p_1\Delta$ and~$p_1 \widehat{\Delta}$ are Ehrhart equivalent.
	
	Finally we consider the geometric mutation of~$t\Delta$ at~$tp_3$. Applying the integral map~$f\defeq\psi_{p_1}^{p_2,p_3}$ defined in Remark~\ref{equimarkov} to~$\Delta$ yields the Markov triangle~$f(\Delta)$ repre\-senting~$(p_1,p_2,p_3)$ in standard~$p_1$-position with the~$p_3$-index lying on the~$x$-axis. We repeat the argument of the mutation at~$tp_2$ after swapping~$p_2$ and~$p_3$ and replacing~$q_2$ by~$q_3\defeq \frac{-q_1p_3-3p_2}{p_1}$, \mbox{and conclude that~$M_{tp_3,tf(\Delta)}(\Z^2)=\Z^2$ when~$t$ divides~$p_1$.}
\end{proof}

\begin{lemma}\label{ancestry}
	Let~$m_n$ be a Markov number occurring as the largest element in a triple~$\mathfrak{m}$ of the~$n$-th generation of the Markov tree, and let~$\Delta$ be a Markov triangle re\-presenting~$\mathfrak{m}$ in standard~$m_n$-position. There exist~$n-1$ Markov numbers~$m_1,\cdots,m_{n-1}$ less than~$m_n$ such that the~$m_1\cdots m_n$-dilate of~$\Delta$ is Ehrhart equivalent to the~$m_1\cdots m_{n}$-dilation of a triangle~$\Delta_{(1,1,1)}$ representing~$(1,1,1)$ in standard~$1$-position. 
\end{lemma}

\begin{proof}
	We prove this by induction on~$n$; the claim holds trivially for~$n=1$. Assume it holds for all positive integers~$k\leq n$, and let~$\mathfrak{m}\eqdef(m_{n-1}',m_n,m_{n-1})$, where~$m_{n-1}'\leq m_{n-1}\leq m_n$. Geometrically mutating~$\Delta$ at the~$m_n$-vertex yields a Markov triangle~$\widehat{\Delta}$ representing the parent~$(m_{n-1}',m_{n-1},\widehat{m_n})$ of~$(m_{n-1}',m_n,m_{n-1})$ in the Markov tree, and the~$[m_{n-1},m_{n-1}']$-edge of~$\widehat{\Delta}$ lies on the~$x$-axis. In particular~$\widehat{\Delta}$ and~$\Delta$ are Ehrhart equivalent because the geometric mutation~$M_{m_n}=M_{(0,0)}$ is integral. We note that the line segment from the~$m_{n-1}$-vertex of~$\widehat{\Delta}$ to the origin is the image of the~$[m_{n-1},m_{n}]$-edge of~$\Delta$ under the geometric mutation, and hence the segments have equal side length, namely~$\frac{m_{n-1}'}{m_{n-1}m_{n}}$. Without loss of generality, we assume that it is the~$m_{n-1}'$-vertex of~$\Delta$ that lies on the~$x$-axis.Then the~$m_{n-1}$-vertex of~$\widehat{\Delta}$ has coordinates~$(\frac{-m_{n-1}'}{m_{n-1}m_{n}},0)$ , from which it follows that translating~$\widehat{\Delta}$ along~$(\frac{m_{n-1}'}{m_{n-1}m_{n}},0)$ yields a Markov triangle~$\Delta'$ representing~$(m_{n-1}',m_{n-1},\widehat{m_n})$ in standard~$m_{n-1}$-position (see Figure~\ref{mutating}). 
	Hence~$m_{n-1}m_{n}\widehat{\Delta}$ is integrally congruent to the~$m_{n-1}m_{n}\Delta'$, so that~$m_{n-1}m_{n}\Delta'$ is Ehrhart equivalent to~$m_{n-1}m_{n}\widehat{\Delta}$ and hence also to~$m_{n-1}m_{n}\Delta$. By the inductive hypothesis applied to~$\Delta'$, there exist~$n-2$ Markov numbers~$m_1,\cdots,m_{n-2}$ smaller than~$m_{n-1}$ such that~$m_1\cdots m_{n-2}m_{n-1}\Delta'$ is Ehrhart equivalent to~$m_1\cdots m_{n-2}m_{n-1}\Delta_{(1,1,1)}$. We conclude that~$m_1\cdots m_{n-2}m_{n-1}m_{n}\Delta$ is Ehrhart equivalent to~$m_1\cdots m_{n-2}m_{n-1}m_{n}\Delta_{(1,1,1)}$.
\end{proof}
	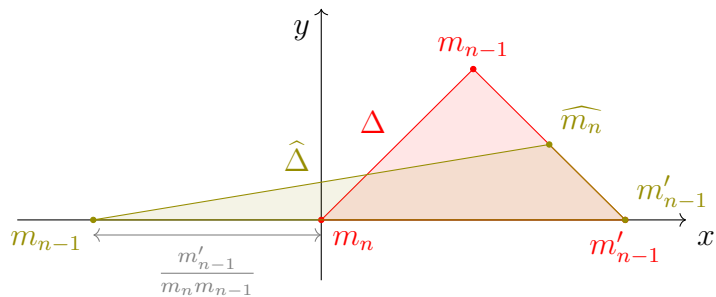
\begin{figure}[h]
		\centering
		\begin{tikzpicture}[scale=2]
			\draw[->] (-2,0) -- (2.4,0) node[below right] {$x$};
			\draw[->] (0,-0.4) -- (0,1.4) node[below left] {$y$};
			\filldraw [red] (0,0) circle [radius=0.5pt] node[anchor=north west] {$m_n$};				
			\filldraw [red] (2,0) circle [radius=0.5pt] node[anchor=north] {$m_{n-1}'$};
			\filldraw [red] (1,1) circle [radius=0.5pt] node[anchor=south]{$m_{n-1}$};
			\filldraw[color=red,fill opacity=0.1](0,0)--(2,0)--(1,1)--cycle;
			\draw[color=red] (0.5,0.5) node[anchor=south east] {$\Delta$};
			
			\filldraw [olive] (-1.5,0) circle [radius=0.5pt] node[anchor=north east] {$m_{n-1}$};
			\filldraw [olive] (2,0) circle [radius=0.5pt] node[anchor=south west] {$m_{n-1}'$};
			\filldraw [olive] (1.5,0.5) circle [radius=0.5pt] node[anchor=south west]{$\widehat{m_n}$};
			\filldraw[color=olive,fill opacity=0.1](-1.5,0)--(2,0)--(1.5,0.5)--cycle;
			\draw[color=olive] (0,0.25) node[anchor=south east] {$\widehat{\Delta}$};
			\draw[color=gray,<->] (-1.5,-0.1)--(0,-0.1);
			\draw[color=gray] (-0.75,-0.1) node[anchor=north] {$\frac{m_{n-1}'}{m_nm_{n-1}}$};
		\end{tikzpicture}
		\caption{Geometric mutation of~$\Delta$ at~$m_n$}
		\label{mutating}
	\end{figure}

\begin{proof}[Proof of Proposition~\ref{periodcollapse}]
	Let~$\Delta$ be a~$(p_1,p_2,p_3)$-triangle in standard~$p_1$-position with some companion number~$q_1\equiv 3p_3p_2^{-2}$~mod~$p_1$ of~$p_1$.
	
	\begin{assertion}\label{assertion_1}
		The period of~$L_\Delta$ is a divisor of~$3p_1$.
	\end{assertion}
	\noindent\textit{Proof of Assertion \ref{assertion_1}.} It suffices to show that the~$3p_1$-dilate of~$\Delta$ is Ehrhart equivalent to the~$3p_1$-dilation of a triangle representing~$(p_1,p_2,p_3)$ in standard~$p_1$-barycentric position, which is pseudo-integral by Proposition~\ref{barycentreposition}. To this end we compute the coordinates of the integral barycentre~$\beta$ of~$\Delta$. By the computation in the proof of Proposition~\ref{pseudomarkov}, the primitive vector~$v_1$ along the integral bisector~$\ell_1$ at the~$p_1$-vertex is a positive multiple of~$(q_1,p_1)$. But~$q_1$ is relatively prime to~$p_1$ by Lemma~\ref{mutationcompanion}, so that~$v_1=(q_1,p_1)$. Thus~$\ell_1$ is defined by the equation~$y=\frac{p_1}{q_1} x$. By Remark~\ref{distancebarycentre}, the affine distance from~$\beta$ to each of the edges of~$\Delta$ is preserved by geometric mutations, and is therefore equal to~$1/3$. Hence the affine distance from~$\beta$ to the~$[p_1,p_2]$-edge with primitive vector~$(1,0)$ is
	\[\frac{1}{3}=|\det\big((1,0),\beta\big)|=\beta_y.\]
	Since~$\beta$ lies on~$\ell_1$, we conclude that~$\beta_x=\frac{q_1}{3p_1}$. Thus the~$t$-dilate of~$\Delta$ is the translation of the~$t$-dilation of a~$(p_1,p_2,p_3)$-triangle~$\Delta'$ in standard~$p_1$-barycentric position along the vector~$(\frac{q_1t}{3p_1},\frac{t}{3})$. This translation preserves the lattice whenever~$3p_1$ divides~$t$, so~$t\Delta$ and~$t\Delta'$ have the same number of integer lattice points if~$3p_1\mid t$. In other terms,~$3p_1\Delta$ and~$3p_1\Delta'$ are Ehrhart equivalent.
	
	\begin{assertion}\label{assertion_2}
		The period of~$L_\Delta$ is a divisor of~$kp_1$, where~$k$ is a product of finitely many Markov numbers which are smaller than~$p_1$.
	\end{assertion}
	\noindent\textit{Proof of Assertion \ref{assertion_2}.} Let~$(p_2',p_1,p_3')$ be the root of the bivalent tree~$\Lambda(p_2',p_1,p_3')$ to which~$(p_1,p_2,p_3)$ belongs, and let~$n$ be the generation of~$(p_2',p_1,p_3')$ and~$m$ the generation of~$(p_1,p_2,p_3)$ in the Markov tree. Applying~$m-n$ geometric mutations to~$\Delta$ (corresponding to decreasing mutations of Markov triples) yields a Markov triangle~$\Delta_0$ representing~$(p_2',p_1,p_3')$ in standard~$p_1$-position. By Lemma~\ref{ehrhartmutate},~$p_1\Delta_0$ and~$p_1\Delta$ are Ehrhart equivalent. By Lemma~\ref{ancestry},~$k\Delta_0$ is Ehrhart equivalent to~$k\Delta_{(1,1,1)}$ where~$k$ is the product of~$n-1$ Markov numbers smaller than~$p_1$; it follows that~$kp_1\Delta$ is Ehrhart equivalent to~$kp_1\Delta_{(1,1,1)}$, which is integral. Therefore the period of~$L_\Delta$ divides~$kp_1$.
	
	By Assertions~\ref{assertion_1}~and~\ref{assertion_2}, the period of~$L_\Delta$ is a common divisor of~$3p_1$ and~$kp_1$, and hence divides~$\operatorname{gcd}(3p_1,kp_1)$. Since every Markov number is relatively prime to~$3$ and~$k$ is a product of Markov numbers, it follows that~$\operatorname{gcd}(3p_1,kp_1)=p_1$, as claimed.
\end{proof}

\begin{proof}[Proof of Corollary~\ref{barnous}]
	By Proposition~\ref{periodcollapse}, the period of such a triangle divides~$2$. The claim follows by observing that the triangle~$\Delta$ with vertices~$(0,0)$,~$(1/2,2)$, and~$(1/2,0)$ is not pseudo-integral, as seen by computing~$L_\Delta(t)$ for~$t=1,2,3$. It follows that the period cannot be~$1$ and must therefore be~$2$.
\end{proof}

\subsection{On the Ehrhart theory of irrational Markov triangles}\label{section_24}
We begin by asserting the Hausdorff convergence of the sequences of Markov triangles in standard~$p_1$-position. Let~$e_\infty$ be the~$[b_\infty,c_\infty]$-edge of~$\Delta_\infty^a$.
\begin{proof}[Proof of Proposition~\ref{hausdorffconv}]
	Let~$C\defeq\operatorname{Cone}(v_a^{(1)},v_a^{(2)})$ be the cone generated by the primitive vectors emanating from the~$a$-vertex of~$\Delta_\infty^a$ (or of~$\Delta_n^a$ for any~$n\in\N$, since~$v_a^{(1)},v_a^{(2)}$ are preserved by the geometric mutations involved in the definition of~$\Delta_n^a$).
	For any~$\e>0$, let~$U_\e$ be the~$\e$-neighbourhood of~$e_\infty$ inside~$C$, i.e.\ the set of all points in~$C$ lying at a \textit{Euclidean} distance smaller than~$\e$ from~$e_\infty$. The limits~\eqref{limitlength} mean that the sequence of~$b_n$-vertices of~$\Delta_n^a$ converges to the~$b_\infty$-vertex of~$\Delta^a_\infty$, and similarly the sequence of~$c_n$-vertices of~$\Delta_n^a$ converges to the~$c_\infty$-vertex of~$\Delta^a_\infty$. In particular for our choice of~$\e$ there exists~$N\in\N$ for which the~$b_n$- resp.~$c_n$-vertex lies within a distance of~$\e$ from the~$b_\infty$- resp.~$c_\infty$-vertex and hence the~$[b_n,c_n]$-edge lies entirely inside~$U_\e$ for~$n\geq N$. Since the complements of~$U_\e$ in our triangles agree for~$n$ sufficiently large, namely~$\Delta_\infty^a\setminus U_\e=\Delta_n^a\setminus U_\e$ for~$n\geq N$, it follows that the symmetric difference~$\Delta_n^a\ominus\Delta_\infty^a\subset U_\e$ is entirely contained in~$U_\e$, which is arbitrarily small.
\end{proof}
	
	\begin{figure}[h]
		\centering
		\begin{tikzpicture}[scale=2]
			\filldraw [blue] (0,0) circle [radius=0.5pt] node[anchor=east] {$a$};
			
			\filldraw [blue] (4,0) circle [radius=0.5pt] node[anchor=north west] {$b_\infty$};
			\filldraw [blue] (2,2) circle [radius=0.5pt] node[anchor= east]{$c_\infty$};
			\draw[color=purple, thick](4,0)--(2,2);
			
			\filldraw [blue] (3.9,0) circle [radius=0.5pt] node[anchor=north east] {$b_n$};
			\filldraw [blue] (2.05,2.05) circle [radius=0.5pt] node[anchor=south east]{$c_n$};
			\draw[color=red](3.9,0)--(2.05,2.05);
			
			\filldraw[color=olive, fill opacity=0.5, dashed] (3.8,0) -- (4.2,0) -- (2.1,2.1) -- (1.9,1.9) -- cycle;
			\draw[color=olive] (3.2,1.2) node {$U_\e$};
			
			\filldraw[color=gray, fill opacity=0.05] (4.6,0) -- (0,0) -- (2.3,2.3);
			\draw[color=gray] (0.5,0.25) node {$C$};
			\draw[color=black,->] (0,0)--(1,0) node[anchor=north] {$v_{a,b_\infty}$};
			\draw[color=black,->] (0,0)--(1,1) node[anchor=south east] {$v_{a,c_\infty}$};
		\end{tikzpicture}
		\caption{Representation of a neighbourhood of the symmetric difference $\Delta_n^a\ominus\Delta_\infty^a$}
		\label{ali}
	\end{figure}
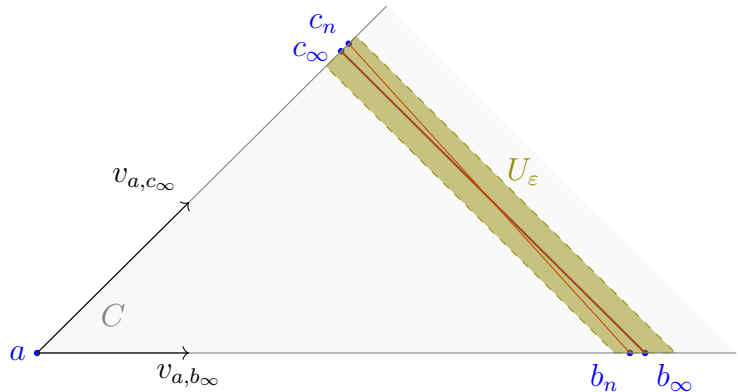\qedhere

\begin{remarkqed}\label{hausdorffnotehrhart}
	Hausdorff convergence does not imply Ehrhart equivalence: given a sequence~$T_n$ of triangles Hausdorff converging to some triangle~$T_\infty$, a lattice point may lie in the exterior of~$T_n$ for all~$n\in\N$ and on the boundary of~$T_\infty$, as shown by the following example. Consider the sequence~$(T_n)$ where~$T_n$ is the triangle with vertices 
	\[\Bigg\{\Big(0,3\frac{\operatorname{Fib}_{2n+1}}{\operatorname{Fib}_{2n+3}}\Big),\Big(3\frac{\operatorname{Fib}_{2n+3}}{\operatorname{Fib}_{2n+1}},0\Big),\Big(3\frac{\operatorname{Fib}_{2n+3}}{\operatorname{Fib}_{2n+1}},3\frac{\operatorname{Fib}_{2n+1}}{\operatorname{Fib}_{2n+3}}\Big)\Bigg\}\]
	which Hausdorff converges to the triangle~$T_\infty$ with vertices
	\[\Big\{(0,\frac{3}{\tau^2}),(3\tau^2,0),(3\tau^2,\frac{3}{\tau^2})\Big\}.\]
	Nevertheless, for all~$n\in\N$,~$T_\infty$ and~$T_n$ are \textit{not} Ehrhart equivalent, since~$T_n$ contains the integer lattice points~$(k,1)$ for~$k=2,\cdots,7$, while~$T_\infty$ contains in addition the point~$(1,1)$.
\end{remarkqed}

Remark~\ref{hausdorffnotehrhart} tells us that to understand the Ehrhart function of~$\Delta_\infty^a$, we must investigate its boundary integer lattice points. The edge~$e_\infty$ spans the affine line 
\[D_\infty: a^2x+\big(a^2\lambda(a)^2-(aq-1)\big)y=a^2\lambda(a).\]
Since~$9a^2-4$ is not a perfect square,~$\sqrt{9-4/a^2}$ is irrational. Therefore the slope of~$D_\infty$ is irrational and hence~$D_\infty$ contains exactly one rational point. We identify it in the following

\begin{lemma}\label{latticepointirrational}
	The only rational point of the affine line~$D_\infty$ is~$\beta\defeq(\frac{q}{3a},\frac{1}{3})$, namely the integral barycentre of~$\Delta_n$ for all~$n\in\N$.
\end{lemma}
\begin{proof}
	Direct computation shows that the equation
	\[a^2x+\big(a^2\lambda(a)^2-(aq-1)\big)y=a^2\lambda(a)\]
	is satisfied by~$(\frac{q}{3a},\frac{1}{3})$.
\end{proof}

\begin{lemma}\label{edgeneighbourhood}
	For each~$t\in\N$,~$te_\infty$ has a neighbourhood~$U_t$ with~$U_t\cap\Z^2\setminus\{t\beta\}=\varnothing$.
\end{lemma}
\begin{proof}
	Fix~$t\in\N$. Let~$d_{(t)}:\R^2\rightarrow\R_{\geq 0}$ be the continuous function mapping each point of the plane to its Euclidean distance from~$tD_\infty$. 
	$d_{(t)}$ clearly attains its minimum~$\e\defeq \min\{d_{(t)}(x,y)\mid (x,y)\in C\cap\Z^2\setminus\beta\}$ over the integer lattice points in the cone~$C$ punctured at~$\beta$. It thus follows from Lemma~\ref{latticepointirrational} that~$\e$ is positive. We define~$U_t\defeq\{(x,y)\in C\setminus\beta\mid d_{(t)}(x,y)<\e\}$ as the~$\e$-tubular neighbourhood of~$te_\infty$~in~$C$. By construction~$U_t\setminus \{t\beta\}$ contains no integer lattice point.
\end{proof}

Lemma~\ref{edgeneighbourhood} allows us to identify the integer lattice points lying in~$t\Delta_\infty^a$:

\begin{prop}\label{latticepointlimit}
	For each~$t\in\N$, there is an~$N_t\in\N$ such that~$t\Delta_\infty^a\cap\Z^2=t\Delta_n^a\cap\Z^2$ for all~$n\geq N_t$.
\end{prop}
\begin{proof}
	Fix~$t\in\N$, and consider the neighbourhood~$U_t$ of~$e_\infty$ given by Lemma~\ref{edgeneighbourhood} so that~$U_t\cap\Z^2\setminus\{t\beta\}=\varnothing$. As in the proof of Proposition~\ref{hausdorffconv}, we see that the sequence of~$t$-dilated Markov triangles~$(t\Delta_n^a)$ that this sequence Hausdorff converges to the~$t$-dilated limiting triangle~$t\Delta_\infty^a$, and therefore there exists~$N_t\in\N$ such that the symmetric difference~$t\Delta_n^a\ominus t\Delta_\infty^a$ is entirely contained in~$U_t$. We note that the point~$t\beta$ belongs to both~$t\Delta_n^a$ as its integral barycentre and to~$t\Delta_\infty^a$ by Lemma~\ref{latticepointirrational}.  Thus the complements~$t\Delta_n^a\setminus (U_t\setminus\{t\beta\})$ and~$t\Delta_\infty^a\setminus (U_t\setminus\{t\beta\})$ of the punctured neighbourhood are equal and hence contain the same integer lattice points. On the other hand both intersections~$t\Delta_n^a\cap (U_t\setminus\{t\beta\})$ and~$t\Delta_\infty^a\cap (U_t\setminus\{t\beta\})$ contain no lattice point by definition of~$U_t$.
\end{proof}

\begin{proof}[Proof of Theorem~\ref{limittriangle}]
	By Proposition~\ref{latticepointlimit}, for any~$t\in\N$,~$at\Delta_\infty^a$ and~$at\Delta_{N_t}^a$ have the same integer lattice points for some~$N_t\in\N$, i.e.~$L_{a\Delta_{N_t}^a}(t)=L_{a\Delta_\infty^a}(t)$. Since all~$a$-dilations of the triangles~$\{\Delta_n^a\}_{n\in\N}$ are Ehrhart equivalent by the proof of Proposition~\ref{periodcollapse},~$L_{a\Delta_n^a}(t)$ does not depend on~$n$, so that~$L_{a\Delta_\infty^a}(t)=L_{a\Delta_n^a}(t)$ for all~$n\in\N$. Since this holds for all~$t$, we conclude that~$L_{a\Delta_\infty^a}=L_{a\Delta_1^a}$.
\end{proof}
\begin{proof}[Proof of Theorem~\ref{limitingbarycentre}]
	We denote the~$[b_\infty,c_\infty]$-edge of~$\Delta_\infty^\beta$ by~$e_\infty^\beta$, which spans the line~$D_\infty^\beta$. By Lemma~\ref{latticepointirrational}, the only rational point of~$D_\infty^\beta$ is the origin~$\beta=(0,0)$. Arguing similarly as in the proof of Lemma~\ref{edgeneighbourhood}, we construct a neighbourhood~$U_t$ of~$te_\infty^\beta$ whose unique integer lattice point is~$(0,0)$. Again arguing similarly as in the proof of Proposition~\ref{latticepointlimit} implies, for each~$t\in\N$, the existence of~$N_t\in\N$ such that~$t\Delta_\infty^\beta\cap\Z^2=t\Delta_n^\beta\cap\Z^2$ for all~$n\geq N_t$. Finally the proof of Theorem~\ref{limittriangle} directly applies in this case as well and we are done.
\end{proof}
	\newpage
	\printbibliography
\end{document}